\documentclass[11pt, twoside, nosumlimits]{amsart}
\usepackage{texdraw}
\usepackage{amssymb}

\def\??#1\\{\marginpar{{\rightskip0pt plus 1fil\bf#1}}}

\def\Jobname{\jobname}
\chardef\catamp=\the\catcode`\@
\catcode`\@=11
\def\t@xdpsfn #1{%
\global\advance \t@xdnum by 1
\ifnum \t@xdnum<10
\xdef #1{\Jobname-00\the\t@xdnum.eps}
\else
\ifnum \t@xdnum<100
\xdef #1{\Jobname-0\the\t@xdnum.eps}
\else
\xdef #1{\Jobname-\the\t@xdnum.eps}
\fi
\fi
}
\catcode`\@=\catamp

\def\Rad{0.1}
\def\factor{0.6}
\def\unit{cm}
\newdimen\Radius
\Radius=\Rad\unit
\Radius=\factor\Radius
\everytexdraw{%
\drawdim pt
\linewd 0.2
\drawdim{\unit}
\setunitscale{\factor}
}
{\catcode`\p=12\catcode`\t=12\gdef\sri#1pt{#1}}
\newbox\txt
\newdimen\Hoehe
\def\Zeichen#1{\htext{\BOX{#1}}}
\def\BOX#1{\setbox\txt=\vbox{\baselineskip6pt%
\halign{$\scriptstyle##$\hfil\cr#1\cr}}%
\global\edef\B{\expandafter\sri\the\wd\txt}%
\global\Hoehe=\ht\txt
\global\advance\Hoehe\dp\txt
\global\edef\H{\expandafter\sri\the\Hoehe}%
\copy\txt}
\chardef\catamp=\the\catcode`\@
\catcode`\@=11
\def\realadd #1#2#3{\dimen0=#1pt
\dimen2=#2pt
\advance \dimen0 by \dimen2
\edef #3{\expandafter\c@lean\the\dimen0}}
\catcode`\@=\catamp
\let\Move=\move
\def\move(#1 #2){\realadd{#1}{1}\Xcoord\realadd{#2}{1}\Ycoord\Move({\Xcoord} {\Ycoord})}
\def\punkt(#1,#2){\move(#1 #2)\fcir f:0 r:\Rad}
\def\kreis(#1,#2){\move(#1 #2)\lcir r:\RRad}
\def\Kreis(#1,#2){\move(#1 #2)\lcir r:\Rad}
\def\strichRO(#1,#2){\move(#1 #2)\rlvec(1 1)}
\def\strichRU(#1,#2){\move(#1 #2)\rlvec(.93 -.93)}
\def\strichU(#1,#2){\move(#1 #2)\rlvec(0 -.9)}
\def\labelO(#1,#2)[#3]#4{%
\punkt(#1,#2)
\rmove(#3 0)
\textref h:C v:B
\Zeichen{#4\vrule depth 5pt width 0pt}
{\drawdim pt
\setunitscale 1
\rmove(0 {\H})
\setunitscale 0.5
\rmove({\B} 0)
\setunitscale 1
\rmove({-\B} 0)
}}
\def\labelU(#1,#2)[#3]#4{%
\punkt(#1,#2)
\rmove(#3 -0.15)
\textref h:C v:T
\Zeichen{#4}
{\drawdim pt
\setunitscale 1
\rmove(0 {-\H})
\setunitscale 0.5
\rmove({\B} 0)
\setunitscale 1
\rmove({-\B} 0)
}\ifx0#2\realmult{0.15}{\factor}\fff
\advance\Hoehe by \fff\unit
\ifdim\Hoehe > \maxT \global\maxT=\Hoehe\fi\fi}
\def\labelR(#1,#2)[#3]#4{%
\punkt(#1,#2)
\rmove(.2 #3)
\textref h:L v:C
\Zeichen{#4}
{\drawdim pt
\setunitscale 1
\rmove({\B} 0)
\setunitscale 0.5
\rmove(0 {\H})
\setunitscale 1
\rmove(0 {-\H})
}\ifx0#2\Hoehe=0.5\Hoehe
\realmult{#3}{\factor}\fff
\advance\Hoehe by -\fff\unit
\ifdim\Hoehe > \maxT \global\maxT=\Hoehe\fi\fi}
\def\labelL(#1,#2)[#3]#4{%
\punkt(#1,#2)
\rmove(-.2 #3)
\textref h:R v:C
\Zeichen{#4}
{\drawdim pt
\setunitscale 1
\rmove({-\B} 0)
\setunitscale 0.5
\rmove(0 {\H})
\setunitscale 1
\rmove(0 {-\H})
}\ifx0#2\Hoehe=0.5\Hoehe
\realmult{#3}{\factor}\fff
\advance\Hoehe by -\fff\unit
\ifdim\Hoehe > \maxT \global\maxT=\Hoehe\fi\fi}
\def\KlabelU(#1,#2)[#3]#4{%
\Kreis(#1,#2)
\rmove(#3 -0.15)
\textref h:C v:T
\Zeichen{#4}
{\drawdim pt
\setunitscale 1
\rmove(0 {-\H})
\setunitscale 0.5
\rmove({\B} 0)
\setunitscale 1
\rmove({-\B} 0)
}\ifx0#2\realmult{0.15}{\factor}\fff
\advance\Hoehe by \fff\unit
\ifdim\Hoehe > \maxT \global\maxT=\Hoehe\fi\fi}
\def\KlabelR(#1,#2)[#3]#4{%
\Kreis(#1,#2)
\rmove(.2 #3)
\textref h:L v:C
\Zeichen{#4}
{\drawdim pt
\setunitscale 1
\rmove({\B} 0)
\setunitscale 0.5
\rmove(0 {\H})
\setunitscale 1
\rmove(0 {-\H})
}\ifx0#2\Hoehe=0.5\Hoehe
\realmult{#3}{\factor}\fff
\advance\Hoehe by -\fff\unit
\ifdim\Hoehe > \maxT \global\maxT=\Hoehe\fi\fi}
\def\labelcond(#1,#2)[#3]#4{%
\move(#1 #2)
\rmove(#3 -0.15)
\textref h:C v:T
\Zeichen{#4}
{\drawdim pt
\setunitscale 1
\rmove(0 {-\H})
\setunitscale 0.5
\rmove({\B} 0)
\setunitscale 1
\rmove({-\B} 0)
}\ifx0#2\realmult{0.15}{\factor}\fff
\advance\Hoehe by \fff\unit
\ifdim\Hoehe > \maxT \global\maxT=\Hoehe\fi\fi}
\def\Vstern(#1,#2){%
\move(#1 #2)
\rmove(-0.1 0.5)
\textref h:R v:C
\Zeichen{\textstyle*}
{\drawdim pt
\setunitscale 1
\rmove({-\B} 0)
\setunitscale 0.5
\rmove(0 {\H})
\setunitscale 1
\rmove(0 {-\H})
}}
\def\Lstern(#1,#2){%
\move(#1 #2)
\rmove(-0.8 0.4)
\textref h:C v:C
\Zeichen{\textstyle*}}
\def\Mstern(#1,#2){%
\move(#1 #2)
\rmove(0 -0.4)
\textref h:C v:C
\Zeichen{\textstyle*}}
\newdimen\maxT
\newbox\Zbox
\long\def\bdiadraw#1\ediadraw{%
\maxT=0pt\setbox\Zbox=\hbox{\btexdraw #1\etexdraw}%
\ifdim\maxT>\Radius\advance\maxT by -\Radius\fi
\lower\maxT\box\Zbox}
\def\labelOO(#1,#2)[#3]#4{%
\move(#1 #2)
\rmove(#3 0)
\textref h:C v:B
\Zeichen{#4\vrule depth 5pt width 0pt}
{\drawdim pt
\setunitscale 1
\rmove(0 {\H})
\setunitscale 0.5
\rmove({\B} 0)
\setunitscale 1
\rmove({-\B} 0)
}}
\def\labelUU(#1,#2)[#3]#4{%
\move(#1 #2)
\rmove(#3 -0.15)
\textref h:C v:T
\Zeichen{#4}
{\drawdim pt
\setunitscale 1
\rmove(0 {-\H})
\setunitscale 0.5
\rmove({\B} 0)
\setunitscale 1
\rmove({-\B} 0)
}\ifx0#2\realmult{0.15}{\factor}\fff
\advance\Hoehe by \fff\unit
\ifdim\Hoehe > \maxT \global\maxT=\Hoehe\fi\fi}
\def\labelUu(#1,#2)[#3]#4{%
\punkt(#1,#2)
\rmove(#3 -0.15)
\textref h:C v:T
\Zeichen{#4}
{\drawdim pt
\setunitscale 1
\rmove(0 {\H})
\setunitscale 0.5
\rmove({\B} 0)
\setunitscale 1
\rmove({-\B} 0)
}}
\def\labelRR(#1,#2)[#3]#4{%
\move(#1 #2)
\rmove(.2 #3)
\textref h:L v:C
\Zeichen{#4}}
\def\labelLL(#1,#2)[#3]#4{%
\move(#1 #2)
\rmove(-.2 #3)
\textref h:R v:C
\Zeichen{#4}}

\catcode`@=11
\newcommand{\REF}[1]{\expandafter\ifx\csname r@#1\endcsname\relax\else{(\it\romannumeral\ref{#1})}\fi}
\renewcommand{\phi}{\varphi}
\renewcommand{\epsilon}{\varepsilon}
\renewcommand{\theta}{\vartheta}
\renewcommand{\rho}{\varrho}
\newcommand{\inn}{\subseteq}
\newcommand{\isom}{\simeq}

\newcommand{\auf}{\hbox{$\rightarrow\hskip-8pt\rightarrow$}}
\newcommand{\into}{\hookrightarrow}
\newcommand{\pfeil}{\rightarrow}

\def\|#1|{\mathop{\rm#1}\nolimits}
\def\_{\penalty10000-\hskip0pt}
\def\lijst#1{\vbox{\baselineskip8pt\halign{\hfill$\scriptstyle##$\cr #1}}}
\def\*{\llap{$*$}}
\newcommand{\mf}{\mathfrak}
\newcommand{\fa}{\mf{a}}
\newcommand{\fb}{\mf{b}}
\newcommand{\fc}{\mf{c}}
\newcommand{\fd}{\mf{d}}

\newcommand{\fg}{\mf{g}}
\newcommand{\fh}{\mf{h}}

\newcommand{\fk}{\mf{k}}
\newcommand{\fl}{\mf{l}}

\newcommand{\fs}{\mf{s}}
\newcommand{\ft}{\mf{t}}

\newcommand{\fz}{\mf{z}}

\newcommand{\gl}{\mf{gl}}
\renewcommand{\sl}{\mf{sl}}
\renewcommand{\sp}{\mf{sp}}
\newcommand{\so}{\mf{so}}
\newcommand{\CC}{{\mathbb C}}

\newcommand{\msf}{\mathsf}
\newcommand{\ssA}{\msf{A}}
\newcommand{\ssB}{\msf{B}}
\newcommand{\ssC}{\msf{C}}
\newcommand{\ssD}{\msf{D}}
\newcommand{\ssE}{\msf{E}}
\newcommand{\ssF}{\msf{F}}
\newcommand{\ssG}{\msf{G}}
\theoremstyle{plain}
\newtheorem{theorem}{Theorem}[section]
\newtheorem{lemma}[theorem]{Lemma}
\newtheorem{prop}[theorem]{Proposition}
\newtheorem{cor}[theorem]{Corollary}
\theoremstyle{definition}
\newtheorem{defn}[theorem]{Definition}
\newtheorem{remark}[theorem]{Remark}
\newtheorem{remarks}[theorem]{Remarks}
\newtheorem{example}[theorem]{Example}

\begin{document}
\pagestyle{myheadings}

\title{Classification of smooth affine spherical varieties}

\author{Friedrich Knop}
\address{Friedrich Knop\\
Department\ of Mathematics\\
Rutgers University\\
110 Frelinghuysen Road\\
Piscataway NJ 08854-8019\\
USA}
\email{knop@math.rutgers.edu}

\author{Bart Van Steirteghem} 
\address{Bart Van Steirteghem\\
Departamento de Matem\'atica\\
Instituto Superior T\'ecnico\\
1049-001 Lisboa\\
Portugal}
\email{bvans@math.ist.utl.pt} 


\begin{abstract}
  Let $G$ be a complex reductive group. A normal $G$\_variety $X$ is called
  spherical if a Borel subgroup of $G$ has a dense orbit in $X$. Of
  particular interest are spherical varieties which are smooth and
  affine since they form local models for multiplicity free
  Hamiltonian $K$\_manifolds, $K$ a maximal compact subgroup of
  $G$. In this paper, we classify all smooth affine spherical
  varieties up to coverings, central tori, and
  $\CC^{\times}$-fibrations.
\end{abstract}

\maketitle

\section{Introduction}\label{sec:intro}

Let $G$ be a connected reductive group (over $\mathbb{C}$). A
normal $G$\_variety $X$ is called {\em spherical} if a Borel
subgroup $B$ of $G$ has a dense orbit in it. When $X$ is affine, this
is equivalent to every simple $G$\_module appearing at most once in
$\CC[X]$ (\cite{VinKim}). In this article, we classify the smooth
affine spherical varieties, up to coverings, central tori and
$\CC^{\times}$-fibrations (see Tables~\ref{sphermod}
through~\ref{BrionTab}).

Our primary motivation are applications to Hamiltonian $K$\_manifolds
where $K\subseteq G$ is a maximal compact subgroup. These are
symplectic $K$\_manifolds which are equipped with a moment map
$m:M\rightarrow\fk^*$. In \cite{Sja}, Sjamaar has shown that locally,
a Hamiltonian $K$\_manifold is isomorphic to a smooth affine
$G$\_variety. A Hamiltonian $K$\_manifold is called {\em multiplicity
  free} if all symplectic reductions are zero\_dimensional.
Multiplicity free Hamiltonian manifolds have spherical varieties as
local models (Brion, \cite{Brion}). Therefore, our classification also
yields a description of the local structure of multiplicity free
Hamiltonian manifolds.

In \cite{Del}, Delzant conjectured, based on results for
commutative groups and for groups of rank 2, that a compact multiplicity 
free
Hamiltonian $K$\_manifold $M$ is completely determined by its generic
isotropy group and the image $m(M)$ of the moment map. The first named
author was able to reduce this conjecture to a statement about smooth
affine spherical varieties: such a variety is completely determined by
the set of highest weights occurring in its ring of regular functions
(Knop conjecture).

The present classification should constitute a major step toward
verifying Knop's conjecture. Unfortunately, difficulties involving
$\CC^\times$\_factors (see below) prevent us from immediately deducing
the conjecture. The second step will be the subject of a forthcoming
paper.

The starting point of our classification is a theorem of Luna,
\cite{LunaSl}, which implies that a general smooth affine spherical
$G$\_variety is a vector bundle $G\times^HV$ where $H\subseteq G$
is a reductive subgroup and $V$ is an $H$\_module.

In principle, it would suffice to classify the triples $(G,H,V)$ but
certain difficulties make the task infeasible. The problems come from
$\CC^\times$\_factors in $G$ or $H$ and from $H$ being non\_connected
(see Examples~\ref{Ex1} through~\ref{Ex2} for details). Therefore, in
this paper, we just determine all triples $(\fg',\fh', V)$ where $G$
and $H$ are replaced by the semisimple part of their Lie algebras. In
the process, we have lost some information. The original triple
$(G,H,V)$ can be recovered with additional combinatorial data but this
point is not addressed in the present paper.

In our classification we build upon existing classifications of two
``extreme'' cases. The first is the case $V=0$, i.e., $X=G/H$ is a
homogeneous variety. Kr\"amer~\cite{Kra} found a complete list of
those when $G$ is simple. His classification was later extended by
Brion~\cite{Br} and Mikityuk~\cite{Mik} to arbitrary $G$. The other
extreme case is when $H=G$, i.e., when $X=V$ is a $G$-module. We call
these {\it spherical $G$\_modules}\footnote{This name collides with
  the notion of a ``spherical representation'', i.e., a representation
  with a $K$\_fixed vector. These will never occur in this paper.}.
Another name for them is ``multiplicity free spaces.''  Kac,
\cite{Kac}, classified irreducible spherical $G$\_modules and Brion,
\cite{Br0}, the reducible spherical modules of simple groups.
Later Benson\_Ratcliff~\cite{BeRa} and Leahy~\cite{Le} independently
completed the classification of spherical modules.

Now, for $G\times^HV$ to be spherical it is necessary that $G/H$ be
spherical for $G$ and that $V$ be spherical for $H$. The converse is
not true. Therefore, we first derive a manageable criterion for
$G\times^HV$ to be spherical. More precisely, for the pair $(G,H)$ we
define a certain subgroup $L$ of $H$, the principal subgroup, with the
property that $G\times^HV$ is spherical if and only if $G/H$ is
spherical and $V$ is spherical for $L$ (Theorem~\ref{crit}).  We
compute the principal subgroup for every variety in the
Kr\"amer\_Brion\_Mikityuk list. Combined with a couple of other useful
``tricks'' the classification turns out to be not too difficult.

While this work was in progress, Camus, \cite{Cam}, independently
classified all smooth affine spherical varieties but only for groups
of type $\ssA$, i.e., when $G$ is locally isomorphic to a product of
$\CC^\times$- and $SL(n)$\_factors. His method is similar to ours
except that we also use the first author's theory of actions on cotangent
bundles, \cite{WM, AsBh}. This simplified matters a lot, see e.g. the
computation of the principal subalgebra, section~\ref{sec:tableproof},
or the reflection ``trick'' of section~\ref{trick}.

Camus, on the other hand, overcame the above mentioned problems with
$\CC^\times$\_factors by using a powerful theory of Luna \cite{Luna}
and thereby deduced Knop's conjecture for groups of type
$\ssA$. Unfortunately, Luna's theory is only verified for groups of
type $\ssA$ and certain other groups and Camus' method does
not immediately carry over to the general case.

\subsection*{Acknowledgments}
We thank the referees for pointing out errors in the proofs of
Lemma~\ref{cutting} and~\ref{so5so6} in a previous version of this paper.
Bart Van Steirteghem was supported by the Fund for Scientific Research
-- Flanders (``aspirant''), by the G.R.F.\ program of the National
Science Foundation (US) and by project FCT/POCTI/FEDER
(Portugal). This was part of his PhD research conducted at Columbia
University under the supervision of Friedrich Knop.

\subsection*{Notation} All varieties we  consider are defined over $\CC$. $G$
will always represent a connected reductive algebraic group. The Lie
algebra of a group is always denoted by the corresponding fraktur
letter.  If $\fg$ is a Lie algebra, we put $\fg':=[\fg, \fg]$.
Fundamental weights $\omega_i$ are numbered as by Bourbaki~\cite{Bou}.
A representation of a semisimple Lie algebra is described by its
highest weight written multiplicatively. For example,
$(\sl(n)+\sp(2k), \omega_1^2+\omega_1\omega_1')$ represents the
$\sl(n)+\sp(2k)$-module $S^2\CC^n \oplus (\CC^n \otimes
\CC^{2k})$. We will sometimes
denote a $\fg$-module $V$ by $(\fg,V)$ to stress which Lie algebra it
is a representation of. To distinguish the $n$\_dimensional abelian
Lie algebra from an $n$\_dimensional representation we denote it by
$\ft^n$. If $V$ is a vector space, then $V^*$ is its dual.

\section{Primitive spherical triples} \label{sec:primitive}

Our starting point is the following well-known application of Luna's Slice
Theorem~\cite[Corollaire 2, p.98]{LunaSl}:

\begin{theorem}\label{lem:hfp}
Let $X$ be a smooth affine $G$\_variety with $\CC[X]^G=\CC$. Then
$X\isom G\times^HV$ where $H$ is a reductive subgroup of $G$ and $V$
is an $H$-module.
\end{theorem}

This applies to spherical varieties in the following way:

\begin{cor} Let $X$ be a smooth affine spherical
  $G$-variety. Then $X\isom G\times^HV$ where $H$ is a reductive
  subgroup of $G$ such that $G/H$ is spherical and $V$ is a spherical
  $H$\_module.
\end{cor}

\begin{proof} Let $B\subseteq G$ be a Borel subgroup. Then $B$ has an
  open orbit in $X$ which implies $\CC[X]^G=\CC$. Moreover, this open
orbit projects to an open orbit in $G/H$ showing that $G/H$ is
spherical. Replacing $B$ by a conjugate if necessary we may assume
that $eH\in G/H$ is in this open orbit. Then $(B\cap H)^\circ$ and, a
fortiori, a Borel subgroup of $H$ has an
open orbit in $V$.
\end{proof}

In the following, we determine (more or less) all triples $(G,H,V)$
such that $G\times^HV$ is spherical. The spherical homogeneous spaces
$G/H$ correspond to the triples $(G,H,0)$ and have been classified by
Kr\"amer~\cite{Kra}, Brion~\cite{Br}, and Mikityuk~\cite{Mik}, while
the spherical $H$\_modules correspond to $(H,H,V)$ and are known thanks
to Kac~\cite{Kac}, Benson-Ratcliff~\cite{BeRa}, and Leahy~\cite{Le}.
The main task of the present paper is therefore to determine all
``mixed'' cases with $H\subsetneq G$ and $V\ne0$.

Next, we present some simple reduction steps.  Given any
two triples $(G_1, H_1, V_1)$ and $(G_2, H_2, V_2)$, we can form their
product $(G_1\times G_2, H_1 \times H_2, V_1 \oplus V_2)$ to obtain a
third. Our classification is therefore one of `indecomposable'
triples. As the following examples show, three other, more subtle
phenomena have to be controlled in order to make the problem
manageable.

\begin{example}[$H$ disconnected]\label{Ex1} 
  Consider the group $G(n):=SL(n)$ and its subgroup $H(n):=SL(n-1)$.
  If $n\ge3$, then the homogeneous space $G(n)/H(n)$ is spherical.
  Now put $G:=G(n_1)\times\ldots\times G(n_s)$ with integers
  $n_1\ge\ldots\ge n_s\ge3$ and $\bar H:=H(n_1)\times\ldots\times
  H(n_s)$.  Then $G/\bar H$ is also a spherical variety which is
  clearly highly decomposable. Now let $N\subseteq G$ be the
  normalizer of $\bar H$ in $G$ (it is isomorphic to
  $GL(n_1-1)\times\ldots\times GL(n_s-1)$). Then $G/H$ is
  spherical too, where $H$ is any subgroup with $\bar H\subseteq H\subseteq
  N$. These subgroups are in one-to-one correspondence to subgroups
  $A$ of $N/\bar H\isom(\CC^\times)^s$. Now choose for $A$ a ``very
  diagonal'' finite subgroup, e.g., the group of $d$\_th roots of unity, $d\ge2$, embedded
  diagonally into $(\CC^\times)^s$. Then the corresponding $H$ is
  disconnected with $H^\circ=\bar H$ making $G/H$ is indecomposable.
\end{example}

Using that $G/H$ is spherical if and only if $G/H^{\circ}$ is we
bypass this problem by considering triples $(\fg, \fh, V)$ where we
have replaced the groups by their Lie algebras.

\begin{example}[$\CC^{\times}$ factors in $H$]
  We keep the notation of the previous example. Now we take $A$ to be
  a ``very diagonal'' connected subgroup, e.g., $\CC^\times$ embedded
  diagonally into $(\CC^\times)^s$. Then again $G/\bar H$ is
  indecomposable. 
\end{example}

\begin{example}[$\CC^{\times}$ factors in $G$]
  Let $G_0:=SL(n_1)\times\ldots\times SL(n_s)$ and
  $G_1:=GL(n_1)\times\ldots\times GL(n_s)$ with $n_1\ge\ldots\ge n_s\ge2$
  acting on $V:=\CC^{n_1}\oplus\ldots\oplus\CC^{n_s}$. Then $V$ is a
  highly decomposable spherical module for both $G_0$ and $G_1$. Now
  let $G$ be any intermediate connected group, determined by a
  connected subgroup of $G_1/G_0\isom(\CC^\times)^s$. Then $V$ will
  in general be an indecomposable spherical variety.
\end{example}

The obvious answer to the problems raised by the last two examples
is to associate to the spherical variety $G \times^H V$ the triple
$(\fg', \fh', V)$ where we replaced the Lie algebras by their
commutator subalgebras. It contains slightly less information. Moreover, we
run into another problem:

\begin{example}[Missing \(\CC^{\times}\) factors]\label{Ex2}
  (1) Let $T\isom\CC^\times$ be a maximal torus in $SL(2)$. To the
  spherical variety $SL(2)/T$ corresponds the triple $(\sl(2), 0,0)$,
  because \(\ft'=0\), but \(SL(2)/e\) is obviously not spherical.
  \newline (2) $V=\CC^n$ is spherical as a $\CC^{\times} \times
  SO(n)$\_variety but not as an $SO(n)$\_variety. Nevertheless, the
  associated triple is $(\so(n), \so(n),\omega_1)$ in both cases.
\end{example}

Taking these issues into account, we define the actual objects of our
classification.

\begin{defn} \label{def:triples}

\begin{enumerate}  
\item Let $\fh\inn\fg$ be semisimple Lie algebras and let $V$ be a
  representation of $\fh$. For $\fs$, a Cartan subalgebra of the
  centralizer $\fc_\fg(\fh)$ of $\fh$, put $\bar\fh:=\fh\oplus\fs$, a
  maximal central extension of $\fh$ in $\fg$. Let $\fz$ be a Cartan
  subalgebra of $\gl(V)^{\fh}$ (the centralizer of $\fh$ in $\gl(V)$).
  We call $(\fg,\fh,V)$ a {\em spherical triple\/} if there exists a
  Borel subalgebra $\fb$ of $\fg$ and a vector $v\in V$ such that
\begin{enumerate}
\item $\fb+\bar\fh=\fg$ and
\item $[(\fb\cap\bar\fh)+\fz]v = V$ where $\fs$ acts via any
  homomorphism $\fs\rightarrow\fz$ on $V$.
\end{enumerate}
\item Two triples $(\fg_i,\fh_i,V_i)$, $i=1,2$, are {\em
isomorphic\/} if there exist linear bijections
$\alpha:\fg_1\pfeil\fg_2$ and $\beta:V_1\pfeil V_2$ such that
\begin{enumerate}
\item$\alpha$ is a Lie algebra homomorphism;
\item$\alpha(\fh_1)=\fh_2$;
\item$\beta(\xi v)=\alpha(\xi)\beta(v)$ for all $\xi\in\fh_1$ and
$v\in V_1$.
\end{enumerate}

\item Triples of the form $(\fg_1\oplus\fg_2,\fh_1\oplus\fh_2,V_1\oplus
V_2)$ with $(\fg_i,\fh_i,V_i)\ne(0,0,0)$ are called {\em
decomposable}.

\item Triples of the form $(\fk,\fk,0)$ and $(0,0,V)$ are said to be {\em
trivial}. 

\item A pair $(\fg,\fh)$ of semisimple Lie algebras is called {\em
spherical} if \((\fg,\fh,0)\) is a spherical triple.

\item A spherical triple (or pair) is {\em primitive\/} if it
  is non-trivial and indecomposable.
\end{enumerate}
\end{defn}

\begin{remarks} 
\begin{enumerate}
\item Recall from the introduction that a $\fg$\_module $V$ is called
  spherical if there is a Borel subalgebra $\fb\subseteq\fg$ and $v\in
  V$ such that $\fb v=V$. We also say that $V$ is spherical for
  $\fg$. This condition is stronger than ``the triple $(\fg,\fg,V)$
  is spherical''. More precisely, the latter is equivalent to $V$ being
  spherical for $\fg+\fz$.

\item If $V=V_1 \oplus \ldots \oplus V_s$ is a decomposition of the
  $\fh$\_module $V$ in irreducible representations, then we can take
  $\fz=\ft^s$ with every $\CC$-factor acting as scalars on the
  corresponding factor $V_i$.

\item The definition of a spherical triple
  is independent of the choice of the Cartan subalgebras $\fs$ and
  $\fz$.

\item For a spherical triple it is very rare that $\fs\ne\fc_\fg(\fh)$ or
  $\fz\ne\gl(V)^\fh$ but is does happen,
  e.g., for $(\sl(2),0,0)$ and
  $(\sl(n),\sl(n),2\omega_1)$, respectively.
  
\item If the first component of two isomorphic triples
  $(\fg_i,\fh_i,V_i)$ is the same, say, $\fg_1=\fg_2=\fg$ then it
  should be kept in mind that $\alpha$ might be an outer
  isomorphism. For example, the triples $(\fh,\fh,V)$ and $(\fh,\fh,V^*)$ are
  always isomorphic in our sense even though $V$ and $V^*$ might be
  different as representations of $\fh$. Other examples are the pairs
  $(\so(8),\so(7))$ and $(\so(8),\mf{spin}(7))$ which are isomorphic due to
  triality.
  
\item From the semisimplicity of $\fg$ and $\fh$ (and $V$) it follows
  that the decomposition of a spherical triple $(\fg, \fh, V)$ into
  indecomposable triples is unique up to isomorphism.

\item A trivial summand of the form $(\fk,\fk,0)$ corresponds to a
    factor $K$ of $G$ which acts trivially on $X$. A trivial summand
    $(0,0,\CC)$ arises if $X$ fibers as a line bundle over a smaller
    space.
\end{enumerate}
\end{remarks}

The next theorem justifies our definition:

\begin{theorem} If $G\times^HV$ is a smooth affine spherical
  variety then $(\fg',\fh',V)$ is a spherical triple. Moreover, it
  follows from the classification that every
  spherical triple arises this way.
\end{theorem}

\begin{proof}
Choose a Borel subgroup $B\subseteq G$ and $v\in V$ such that the
$B$\_orbit of $[1,v]\in G\times^HV$ is dense. This is equivalent to
the density of $BH$ in $G$ and of $(B\cap H)v$ in $V$. In terms of Lie
algebras we get $\fb+\fh=\fg$ and $(\fb\cap\fh)v=V$.

Now let $\fg=\fg'\oplus\fa$ and $\fh=\fh'\oplus\fd$. Since the
conditions of Definition~\ref{def:triples} are independent on the
choice of $\fs$ and $\fz$ we may arrange that the image of $\fd$ in
$\fg'$ and $\fg\fl(V)$ is contained in $\fs$ and $\fz$, respectively.
Let, in abuse of notation, $\fb'=\fb\cap\fg'$, a Borel subalgebra of
$\fg'$. Then the sphericality of $(\fg',\fh',V)$ follows from the following two computations:
$$
\fg'\oplus\fa=(\fb'\oplus\fa)+(\fh'+\fd)\subseteq
(\fb'\oplus\fa)+(\fh'+\fs+\fa)=(\fb'+\fh'+\fs)\oplus\fa
$$
and
$$
V=((\fb'\oplus\fa)\cap(\fh'+\fd))v\subseteq(\fb'\cap(\fh'+\fs)\oplus\fz)v.
$$
The last inclusion is seen as follows: let $h+d\in\fh'+\fd$ with
$h+d\in\fb'\oplus\fa$. Then $(h+d)v\in(h+\fz)v=(h+s+\fz)v$ where
$s\in\fs$ is the projection of $d$ into $\fg'$. Then
$h+s\in\fb'\cap(\fh'+\fs)$ and the inclusion follows.

For the second part, let $G'$ be the simply connected group with Lie
algebra $\fg'$, and let $H'\subseteq G'$ be the connected subgroup
with Lie algebra $\fh'$. Furthermore, let $S$ and $Z$ be maximal tori
in the centralizers $C_{G'}(H')$ and $GL(V)^{H'}$, respectively.

Next we need that the action of $\fh'$ on $V$ integrates to an action
of $H'$ on $V$. Since $H'$ is, in general, not simply connected this
is a non\_trivial condition. It follows a posteriori, by inspection
of Table~\ref{newtab2} and the inference rules (Table~\ref{rules}).

Granted this, we define $G:=G'\times S\times Z$ and $H:=H'\times
S\times Z$. Here $S\subseteq H$ embeds diagonally into $G'\times S$
and acts trivially on $V$. We claim that $G\times^HV$ is
spherical. Indeed,
$\fb+\fh=(\fb'\oplus\fs\oplus\fz)+(\fh'\oplus\Delta\fs\oplus\fz)=
(\fb'+\fh'+\fs)\oplus\fs\oplus\fz=\fg$ and
$(\fb\cap\fh)v=((\fb'\cap(\fh'+\fs))\oplus\fz)v=V$, since $\fs$ acts through $\fz$.
\end{proof}

\section{Diagrams}

To get a hold on the combinatorics involved we generalize the notation
of Mikityuk \cite{Mik} and represent a triple by a
three-layered graph $\Gamma$ as follows: Let
$\fg=\fg_1+\ldots+\fg_r$, $\fh=\fh_1+\ldots\fh_s$,
$V=V_1\oplus\ldots\oplus V_t$ be the decompositions into simple
factors. The vertices of the graph are all $\fh_j$, all $V_k$, and
those $\fg_i$ which are not contained in $\fh$. There is an edge
between $\fh_j$ and $\fg_i$ if $\fh_j\into\fg\auf\fg_i$ is non-zero
and an edge between $V_k$ and $\fh_j$ if $V_k$ is a non-trivial
$\fh_j$\_module. There is no edge between $V_k$ and $\fg_i$.

In practice, we write the graph in three rows with the
$\fg_i$\_vertices, $\fh_j$\_vertices, and $V_k$\_vertices in the first,
second, and third row respectively. In principle, all edges should be
labeled to say how $\fh_j$ is embedded into $\fg_i$ or how $\fh_j$
acts on $V_k$. In most cases these embeddings and actions are the
``natural'' ones and we then omit the labels. Note that the $\fg_i$
which are contained in $\fh$ can be recovered from $\Gamma$: they
correspond to those vertices $\fh_j$ which are not connected to any
$\fg_i$.

\begin{example}
The graph
$$
\begin{texdraw}
\strichU(0,2)
\strichRU(0,1)
\strichRO(1,0)
\punkt(1,0)
\labelO(0,2)[0]{\sl(m{+}1)}
\labelL(0,1)[0]{\sl(m)}
\labelO(2,1)[0]{\sl(n)}
\end{texdraw}
$$
represents the triple
$(\sl(m+1)+\sl(n),\sl(m)+\sl(n),\omega_1\omega_1')$.
\end{example}

\section{The classification} \label{sec:fullclassif}

\begin{theorem}

(1) Every primitive spherical triple $(\fg,\fh,V)$ with $V\ne0$ is
    contained in Table~\ref{sphermod} or~\ref{newtab2} or can be
    obtained from an item in these tables by applying the inference
    rules (Table~\ref{rules}), possibly several times.

(2) If $(\fg,\fh,0)$ is a primitive spherical triple then $(\fg,\fh)$
    is contained in either Table~\ref{KraemerTab} or Table~\ref{BrionTab}.

Moreover, every triple in these tables is spherical.
\end{theorem}

Part (2) of the theorem is just added for completeness since it
restates the classifications of Kr\"amer, Brion, and
Mikityuk. Table~\ref{sphermod} restates the classification of spherical
modules by Kac, Benson\_Ratcliff, and Leahy. New are
Tables~\ref{newtab2} and~\ref{rules}.

\medskip\noindent
{\bf Explanation of the inference rules:} They were introduced to keep
the table at a moderate size. Expanding all items using these rules
would result in 30+ more cases. A circled vertex
$\begin{texdraw}
\punkt(0,0)
\kreis(0,0)
\end{texdraw}$
should remind the
reader that an application of an inference rule is possible at that vertex.

\subsection{}The first rule means that every $(\sl(2),\sl(2))$\_summand of
$(\fg,\fh)$ can be replaced by $(\sp(2m+2),\sl(2)+\sp(2m))$ with
$m\ge1$.

\subsection{}The second rule states that if $(\fg,\fh)$ contains the
summand $(\sp(4),\sp(4))$ and, as an $\sp(4)$\_module, $V$ contains
only the trivial or the defining representation then that summand can
be replaced by $(\so(7),\sl(4))$ (with embedding
$\sl(4)\isom\so(6)\into\so(7)$). Moreover, the representation
$(\sp(4),\omega_1)$ can be replaced by either extension
$(\sl(4),\omega_1)$ or $(\sl(4),\omega_3)$. Since there is an element
of $\so(7)$ which induces the outer automorphism of $\sl(4)$ this
choice is only relevant if $\sp(4)$ acts non\_trivially on more than
one component of $V$. There is exactly one such case and it yields
$$
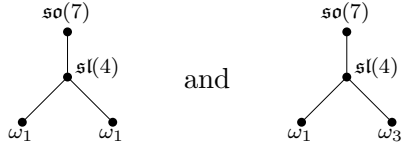
\begin{texdraw}
\strichRO(0,0)
\strichU(1,2)
\strichRU(1,1)
\labelO(1,2)[0]{\so(7)}
\labelR(1,1)[.2]{\sl(4)}
\labelU(0,0)[0]{\omega_1}
\labelU(2,0)[0]{\omega_1}
\end{texdraw}
\qquad\raise20pt\hbox{\rm and}\qquad
\begin{texdraw}
\strichRO(0,0)
\strichU(1,2)
\strichRU(1,1)
\labelO(1,2)[0]{\so(7)}
\labelR(1,1)[.2]{\sl(4)}
\labelU(0,0)[0]{\omega_1}
\labelU(2,0)[0]{\omega_3}
\end{texdraw}
$$

\begin{example} Applying both rules to the same diagram:
\medskip

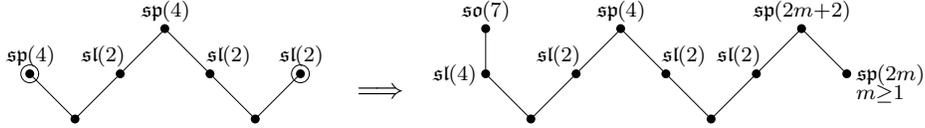
\begin{texdraw}
\strichRU(0,1)
\strichRO(1,0)
\strichRO(2,1)
\strichRU(3,2)
\strichRU(4,1)
\strichRO(5,0)
\kreis(0,1)
\punkt(1,0)
\punkt(5,0)
\kreis(6,1)
\labelO(0,1)[0]{\sp(4)}
\labelO(2,1)[-0.4]{\sl(2)}
\labelO(3,2)[0]{\sp(4)}
\labelO(4,1)[0.4]{\sl(2)}
\labelO(6,1)[0]{\sl(2)}
\end{texdraw}
\quad$\raise10pt\hbox{$\Longrightarrow$}$\quad
\begin{texdraw}
\strichRU(0,1)
\strichRO(1,0)
\strichRO(2,1)
\strichRU(3,2)
\strichRU(4,1)
\strichRO(5,0)
\strichU(0,2)
\labelO(0,2)[0]{\so(7)}
\punkt(1,0)
\punkt(5,0)
\strichRO(6,1)
\strichRU(7,2)
\labelO(7,2)[0]{\sp(2m{+}2)}
\labelR(8,1)[0]{\sp(2m)}
\labelRR(8,0.6)[0]{m\ge1}
\labelL(0,1)[0]{\sl(4)}
\labelO(2,1)[-0.4]{\sl(2)}
\labelO(3,2)[0]{\sp(4)}
\labelO(4,1)[0.4]{\sl(2)}
\labelO(6,1)[-0.4]{\sl(2)}
\end{texdraw}

\end{example}

Another measure we took to cut down on the length of the table was
allowing borderline cases: whenever an $\fh_j$ is $0$, i.e.,
$\fh_j=\sl(n)$, $n\le1$ or $\sp(2n)$, $n\le0$ then that vertex and
all adjacent edges are to be omitted. For that reason,
Table~\ref{BrionTab} contains at first sight fewer entries than
the tables of Brion and Mikityuk. Incidentally, the table of Kr\"amer
contains a couple of redundancies which we removed, e.g.,
$(\so(8),\sp(4)+\sl(2))$ is, using triality, just
$(\so(8),\so(5)+\so(3))$.

\section{Tools}\label{sec:tools}

\subsection{Criterion for sphericality}\label{sec:criterion}
Our main tool is a manageable criterion for a triple to be spherical.
It will replace the `fiber condition' in Definition~\ref{def:triples}
with a sphericality condition on a representation of a {\em
  reductive\/} Lie algebra. We are going to derive the criterion from
\cite{WM,AsBh} but we could also have used \cite{Pan1, Pan2}.

Let $X$ be a $G$\_variety. Then $L\subseteq G$ is called a {\em
generic isotropy subgroup} if $L$ is conjugate to $G_x$ for all $x$ in a
non\_empty open subset of $X$. If one exists, it is unique up to
conjugacy. We are going to use this concept only in case $X$ is a
vector space. In that case, the existence of a generic isotropy
subgroup is guaranteed by theorems of Richardson \cite{Rich} and
Luna~\cite{LunaSl}. The Lie algebra of $L$ is called a {\em generic
isotropy subalgebra}.

\begin{defn} \label{def:princ} Let $(\fg,\fh)$ be a spherical pair and
  $\bar\fh=\fh+\fs\subseteq\fg$ a maximal central extension. Let
$\bar\fl$ be a generic isotropy subalgebra of $\bar\fh$ acting on
$\bar\fh^\perp$, the orthogonal complement of $\bar\fh$ in $\fg$. Then
the image $\fl$ of $\bar\fl$ under the projection
$\fh+\fs\auf\fh$ is called a {\it principal subalgebra\/} of
$(\fg,\fh)$.
\end{defn}

It is known that $\bar\fl$, and hence $\fl$, is always reductive
(\cite[Korollar 8.2]{WM}, \cite{Pan1}). The following lemma
establishes the connection of $\bar\fl$ with Borel subgroups. It is
also contained in \cite{Pan1}.

\begin{lemma}\label{borel}
  Let $H$ be a reductive subgroup of $G$ such that $BH \inn G$ is
  dense. Then $B\cap H$ is a Borel subgroup of a
  generic isotropy group for $H$ acting on $\fh^{\perp}$.
\end{lemma}
\begin{proof} Let $U$ be the open subset of $\fh^{\perp}$ such that
  for every $u\in U$, $H_u$ is $H$\_conjugate to a generic isotropy
  group $L$ and put $x_o:= eH\in G/H$. From \cite[Theorem 3.2]{AsBh}
  we obtain a $B$\_invariant section $\sigma: Bx_o \rightarrow
  T^*_{G/H}=G\times^H \fh^{\perp}$ which intersects $G\times^H U$
  non\_trivially. This implies $v\in U$ where $\sigma(x_o)=[e,v]$.
  Hence
\[
B\cap H = B_{x_o} \inn G_{\sigma(x_o)}=H_v=L.
\]
On the other hand, Corollaries 2.4 and 8.2 in~\cite{WM} imply that
$B\cap H$ is $G$\_conjugate to a Borel subgroup of $L$. Therefore, it
is a Borel subgroup of $L$.
\end{proof}

 The criterion is now:

\begin{theorem}\label{crit}
The triple $(\fg,\fh,V)$ is spherical if and only if $(\fg,\fh)$ is a
spherical pair and $V$ is a spherical $\fl+\fz$\_module. Here
$\fl$ is a principal subalgebra of $(\fg,\fh)$ and $\fz$ is a Cartan
subalgebra of $\gl(V)^\fh$.
\end{theorem}

\begin{proof}
  Assume $(\fg,\fh)$ is a spherical pair. Lemma~\ref{borel} implies
  that $\fb\cap(\fh+\fs)$ is a Borel subalgebra of
  $\bar\fl\subseteq\fh+\fs$. Thus, $(\fg,\fh,V)$ is spherical if
  and only if $V$ is a spherical $\bar\fl+\fz$\_module. The image of $\fb\cap(\fh+\fs)$
  in $\fh$ is a Borel subalgebra of $\fl$. Since $\fs$ acts on $V$
  through $\fz$, we can replace $\bar\fl$ by $\fl$.
\end{proof}

Observe that ``$V$ is a spherical
  $\fl+\fz$\_module'' is a stronger condition than
  ``$(\fl',\fl',V)$ is a spherical triple'', since
  $\fh$\_irreducible modules may not be $\fl$\_irreducible. In other
  words, in the latter statement, $\fz$ could be bigger. To deal with
  this kind of ``$\CC^\times$\_deficiency'' one can use the following

\begin{lemma}\label{deficient}
Let $(\fh,V)$ be a spherical module and $\fz\subseteq\gl(V)^\fh$ a
  Cartan subalgebra. Then there is a (unique) subspace
  $\fc\subseteq\fz$ such that for every subspace $\fz_0\subseteq\fz$:
$$
\hbox{$V$ is spherical for $\fh+\fz_0$}\qquad
\Leftrightarrow\qquad\fz_0+\fc=\fz.
$$
\end{lemma}

For a proof see \cite{KRem}~Thm.~5.1. In the notation of that paper we
have $\fc=(\fa^*\cap\fz^*)^\perp$. For the convenience of the reader
we list $\fc$ for the cases we are going to use in the sequel: we have
$\fc=0$ for
$$
(\fh,V)=(\so(n),\omega_1), (\so(8),\omega_3+\omega_4),
(\sp(2m)+\sl(2)+\sp(2n),\omega_1\omega_1'+\omega_1'\omega_1'').
$$
In other words, for these pairs, $V$ will not remain spherical if
$\fz$ is replaced by any proper subspace $\fz_0$. In the other
extreme, we have $\fc=\fz$ for $(\fh,V)=(\sp(2n),\omega_1)$. Finally,
we have the following intermediate cases:
\begin{align*}
\fc=\CC&(1,1)\hbox{ for }(\fh,V)=(\sl(n),\omega_1+\omega_1), n\ge2.
\\
\fc=\CC&(1,-1)\hbox{ for }(\fh,V)=(\sl(n),\omega_1+\omega_{n-1}), n\ge3.
\\
\fc\subseteq\CC&(1,1)\hbox{ for
}(\fh,V)=(\sl(2)+\sl(n),\omega_1\omega_1'+\omega_1'), n\ge2.
\end{align*}

Note that in the first and the third case the scalar
$\ft^1$\_action on $V$ does {\it not} make $V$ spherical for $\fh\oplus\ft^1$.

\subsection{Calculating principal subalgebras} \label{sec:tableproof}
The next step is to determine the principal subalgebras.

\begin{prop}\label{principal}
  The principal subalgebras of all primitive spherical pairs are
  listed in Tables~\ref{KraemerTab} and~\ref{BrionTab}.
\end{prop}

The verification of the data is quite standard and is left to the
reader. In the calculations, Elashvili's tables in~\cite{ElasCan} and
\cite{ElasStat} are very useful, especially, if $\fs=0$. In all other
cases we have $\fs=\ft^1$ (by inspection). In the tables, these cases are
indicated by a ``$*$''. For them, the following observation proves
useful since Elashvili's tables only contain generic isotropy algebras
for representations of semisimple Lie algebras.

\begin{lemma} Let $(\fg,\fh)$ be a primitive spherical pair with $\fs=\ft^1$, let
  $\fl$ be a principal subalgebra, and let $\fl_*$ be a generic
  isotropy subalgebra of $\fh$ acting on $\fh^\perp$. Then the
  following statements are equivalent:
\begin{enumerate}
\item\label{L1}$\fl\isom\fl_*\oplus\ft^1$.
\item\label{Kcond}$G/H$ is spherical.
\item\label{Scond}Every $\fh$\_invariant on $\bar\fh^\perp$ is
$\bar\fh$\_invariant.
\end{enumerate}
If these conditions do not hold then $\fl=\fl_*$.
\end{lemma}

\begin{proof} First observe that the map $\bar\fl\rightarrow\fl$ is
  bijective. Otherwise $\fs\subseteq\bar\fl$, since $\dim \fs=1$, and $\fs$, being normal
  in $\fh+\fs$, would act
  trivially on $(\fh+\fs)^\perp$. This implies that $\fs$ is in
  the center of $\fg$. Contradiction.

Let $\bar H=S\cdot H$ and consider the $\CC^\times$\_fibration
$X:=G/H\rightarrow \bar X:=G/\bar H$. Let $U\subseteq G$ be a maximal
unipotent subgroup. Then the dimension of a generic $U$\_orbit in $X$ and
$\bar X$ is the same, which implies that $P_u(X)=P_u(\bar X)$ in the
notation of \cite{IB}. Therefore, the sum of rank and complexity is
one smaller for $\bar X$ than for $X$
(\cite[Korollar~2.12]{IB}). Thus, $X$ is spherical if and only if its complexity is zero if and only if
${\rm rk}X={\rm rk}\bar X+1$ if and only if a generic isotropy group
in $T^*_X$ is of codimension one in that of $T^*_{\bar X}$. This
proves the equivalence of \REF{L1} and \REF{Kcond}. The equivalence of
\REF{Kcond} and \REF{Scond} follows from \cite[Satz~7.1]{WM}.
If $X$ is not spherical then ${\rm rk}X={\rm rk}\bar X$ and 
$\fl_*=\bar\fl$.
\end{proof}

\begin{remark} If $\fg$ is simple, then information on~\REF{Kcond} is
  part of Kr\"amer's table. For \REF{Scond}, the tables of
  Schwarz~\cite{Schwarz} are useful.
\end{remark}

\subsection{A reduction lemma}
In certain cases, it is not even necessary to calculate the principal
subalgebra.

\begin{lemma}\label{reductionlemma}
  Let $(\fg,\fh)$ be a spherical pair such that $\bar\fh^\perp=U\oplus
  U^*$ with $U$ an irreducible spherical $\bar\fh$\_module. Let
  $\tilde\fh\subseteq\tilde\fg$ be another pair of semisimple Lie
  algebras and let $V$ be an $\fh+\tilde\fh$\_module. Consider the
  following statements:
\begin{enumerate}
\item The triple $(\fg+\tilde\fg,\fh+\tilde\fh,V)$ is
  spherical.\label{reditem1}
\item The triple $(\fh+\tilde\fg,\fh+\tilde\fh,U\oplus V)$ is
  spherical.\label{reditem2}
\end{enumerate}
Then \REF{reditem1} implies \REF{reditem2}. The converse is true if
the action of $\bar\fh$ on $U$ contains the scalars.
\end{lemma}

\begin{proof} Clearly, both conditions imply that
  $(\tilde\fg,\tilde\fh)$ is a spherical pair. Let $\fl\subseteq\fh$
and $\tilde\fl\subseteq\tilde\fh$ be the principal subalgebras and
$\fz\subseteq\gl(V)^{\fh+\tilde\fh}$ a Cartan subalgebra. Then
\REF{reditem1} holds if and only if $V$ is spherical for
$\fl+\tilde\fl+\fz$. 

Let $\fb_0\subseteq\bar\fh$ be a Borel subalgebra. Since $U$ is
  spherical for $\bar\fh$ there is $u\in U$ with $\fb_0u=U$. Let
  $\fb_1\subseteq\fb_0$ be its isotropy subalgebra. As in
  Lemma~\ref{borel}, $\fb_1$ is a Borel subalgebra of a generic
  isotropy subalgebra $\bar\fl$ for $\bar\fh$ acting on $T^*_U=U\oplus
  U^*=\bar\fh^\perp$. Thus, \REF{reditem1} holds if and only if
  $U\oplus V$ is spherical for $\bar\fh+\tilde\fl+\fz$. This implies
  that $U\oplus V$ is spherical for $\fh+\tilde\fl+\ft^1+\fz$, i.e.,
  \REF{reditem2}. Moreover, the last implication is an equivalence
  if $\ft^1\subseteq\bar\fh$.
\end{proof}

\subsection{The cutting lemma}
Clearly, a triple is trivial if its graph consists of only isolated
vertices which are all in the $\fh$- and $V$\_layers\footnote{Note
  that the non\_trivial triple $(\sl(2),0,0)$ is represented by an
  isolated vertex in the $\fg$\_layer.}. A triple is primitive if and
only if it is non-trivial and its graph is connected.

The pair $(\fg,\fh)$ is called the {\it base} of the triple. Its graph
$\Gamma_b$ is obtained by removing the $V$\_layer from $\Gamma$. A
non\_trivial component of $(\fg,\fh)$ is called a {\it base
component}. It corresponds to a connected component of $\Gamma_b$
which is not an isolated vertex in the $\fh$\_level.

The triple $(\fh,\fh,V)$ is called the {\it fiber} of the triple. Its
graph $\Gamma_f$ is obtained by removing the $\fg$\_layer from
$\Gamma$. A component of the fiber is called a {\it fiber
  component}. It corresponds to a connected component of $\Gamma_f$.

The process of removing the top or bottom layer can be
generalized. Let $v$ be a vertex of $\Gamma$. A {\it cut of $\Gamma$
in $v$} is a graph $\Gamma'$ which has the same vertices as $\Gamma$
except that $v$ is replaced by two vertices $v'$, $v''$ (in the same
layer). Moreover, $\Gamma'$ has the same edges as $\Gamma$ with the
property that each edge adjacent to $v$ is in $\Gamma'$ adjacent to
either $v'$ or $v''$.

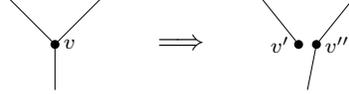
\begin{example}

\

$$
\begin{texdraw}
\labelR(1,1)[0]{v}
\move(1 1)\rlvec(-1 1)
\move(1 1)\rlvec(1 1)
\move(1 1)\rlvec(0 -1)
\end{texdraw}
\qquad\raise 15pt\hbox{$\Longrightarrow$}\qquad
\begin{texdraw}
\labelL(.8,1)[0]{v'}
\labelR(1.2,1)[0]{v''}
\move(.8 1)\rlvec(-.8 1)
\move(1.2 1)\rlvec(.8 1)
\move(1.2 1)\rlvec(-.2 -1)
\end{texdraw}
$$

\end{example}

We need the following

\begin{lemma}\label{Vtensor} Let $U_1\ne0$, $U_2\ne0$, $U_3$ be
  representations of the reductive Lie algebra $\fl$. Assume that
  $(U_1\otimes U_2)\oplus U_3$ is a spherical $\fl+\ft^2$\_module.
  Then $U_1\oplus U_2\oplus U_3$ is a spherical $\fl+\ft^3$\_module.
\end{lemma}

\begin{proof}
  Consider the morphism
$$
\pi:U_1\oplus U_2\oplus U_3\rightarrow (U_1\otimes U_2)\oplus U_3:
(u_1,u_2,u_3)\mapsto(u_1\otimes u_2,u_3).
$$
The image of $\pi$ is a subvariety of a spherical variety, hence
spherical. Moreover, $\pi$ is generically a $\CC^\times$\_fibration.
More precisely, it is the quotient by the $\CC^\times$\_action
$t\cdot(u_1,u_2,u_3)=(tu_1,t^{-1}u_2,u_3)$, of which the infinitesimal action
is contained in that of $\ft^3\subseteq\fl+\ft^3$ on $U_1\oplus U_2\oplus U_3$. This implies that
$U_1\oplus U_2\oplus U_3$ is spherical as well.
\end{proof}

Now we can prove the cutting lemma:

\begin{lemma}
\label{cutting}
  Let $\Gamma$ be the graph of a spherical triple and let $\Gamma'$ be
  the graph obtained by cutting $\Gamma$ in a vertex in the
  $\fh$\_layer or in the $V$\_layer. Then $\Gamma'$ is also the graph of
  a spherical triple.
\end{lemma}

\begin{proof}
1. {\it Cut in the $V$\_layer.} Let $v$ be a vertex corresponding to
   an irreducible component $V_k$ of $(\fh,V)$. The adjacent edges
   correspond to the simple factors of $\fh$ acting non\_trivially on
   $V_k$. Thus there is a decomposition $V_k=U_1\otimes U_2$ such that
   a simple factor $h_j$ which is attached to $v'$ or $v''$ (in $\Gamma'$) acts only on
   $U_1$ or $U_2$, respectively. Write $V=V_k\oplus U_3$. Then the
   process of cutting amounts to replacing $V$ by $U_1\oplus U_2\oplus
   U_3$. The assertion follows from Lemma~\ref{Vtensor}.
   
   2. {\it Cut in the $\fh$\_layer.} Let $v$ be a vertex corresponding
   to a simple factor $\fh_j$ of $\fh$. If both cut vertices $v'$ and
   $v''$ are connetced to the $\fg$\_layer then $\Gamma'$ corresponds
   to the triple $(\fg,\fh\oplus\fh_j,V)$. Let $\overline\fl$ be a
   principal subalgebra of $(\fg,\fh\oplus\fh_j)$.  Then
   $\fl\subseteq\overline\fl$ implies that the image of $\overline\fl$
   in $\gl(V)$ contains the image of $\fl$, which proves the
   assertion. If one of the verices $v'$, $v''$ is not connected to
   the $\fg$\_layer then $\Gamma'$ represents the triple
   $(\fg\oplus\fh_j,\fh\oplus\fh_j,V)$. Its principal subalgebra is
   $\fl\oplus\fh_j$ and we argue as before.
\end{proof}

An immediate consequence is the erasing lemma:

\begin{cor}
\label{erase}
Let $(\fg,\fh,V)$ be a spherical triple with graph $\Gamma$. Let
$\tilde\Gamma$ be the graph obtained from $\Gamma$ by erasing any
number of edges between the $\fh$- and the $V$\_layer. Then the
triple $(\fg,\fh,\tilde V)$ corresponding to
$\tilde\Gamma$ is spherical, as well.
\end{cor}

\begin{proof}
Indeed, we can cut out every edge.
\end{proof}

\section{The base components}\label{subsec:comp}

In this section we start the classification by ruling out most of the
primitive spherical pairs from being base components of a primitive triple with nonzero fiber.

\begin{prop}
\label{ActualComp}
Every primitive spherical pair which is a base component of a primitive
spherical triple $(\fg,\fh,V)$ with $V\ne0$ is contained in
Table~\ref{PairsForTriples}. The top two layers of the graph represent
the component $(\fg_0,\fh_0)$ while the third layer indicates its
principal subalgebra $\fl$ and its embedding into $\fh_0$.
\end{prop}

The rest of this subsection is devoted to the proof of this statement.
Let $(\fg_0,\fh_0)$ be as in the proposition. Then the erasing
Lemma (Corollary~\ref{erase}) implies that there is a primitive spherical triple
of the form $(\fg,\fh,V)$ with $\fg=\fg_0$, $\fh=\fh_0$ and where $V$
is irreducible and exactly one simple factor $\fh_j$ of $\fh$ acts
non\_trivially on $V$.

Some simple algebras have more than one irreducible spherical
module. The following lemma reduces the number of cases to check:
 
\begin{lemma} \label{AAAA} Let $(\fg,\fh,V)$ be a primitive spherical
  triple such that $V$ is irreducible and exaclty one simple factor,
  say $\fh_j$, acts non\_trivially on $V$.
\begin{enumerate}
\item \label{pleth} Assume $\fh_j=\sl(m)$, $m\ge2$, and $V=\omega_1^2$,
  $V=\omega_{m-1}^2$, or $V=\omega_{m-1}$. Then also
  $(\fg,\fh,\omega_1)$ is spherical.
\item \label{sp4} Assume $\fh_j=\sp(4)$ and $V=\omega_2$. Then also
  $(\fg,\fh,\omega_1)$ is spherical.
\end{enumerate}
\end{lemma}

\begin{proof}
Let $\fl$ be the principal subalgebra of $(\fg,\fh)$ and
$\tilde\fl=\fl\oplus\ft^1$. 

\noindent\REF{pleth} It is well known that a module is
  spherical if and only its dual is. Thus we may assume
  $V=\omega_1^2$. Let $\fb$ be a Borel subalgebra of
$\gl(m)$. Since $\dim\omega_1^2=\dim\fb$, the image of $\tilde\fl$ in
$\gl(m)$ contains $\fb$, hence equals $\gl(m)$.

\noindent\REF{sp4} Let $\fl_0$ be the image of $\fl$ in
$\fh_j=\sp(4)$. The Borel subalgebra of $\fl_0$ must have dimension at
least $\|dim|V-1=4$. This leaves only the possibilities $\fl_0=\sp(4)$
or $\fl_0=\sl(2)+\sl(2)$. But in the latter case $V$ is not
spherical for $\fl_0\oplus\ft^1$.
\end{proof}

\begin{proof}[Proof of Proposition~\ref{ActualComp}]
First, the items of Tables~\ref{KraemerTab}
and~\ref{BrionTab} marked by a ``$\bullet$'' are members of
Table~\ref{PairsForTriples}. We examine the others.

Let $\fl$ be a principal subalgebra of $(\fg, \fh)$ and
$\tilde\fl=\fl\oplus\ft^1$. According to the list of irreducible
spherical modules (in Table~\ref{sphermod}) and Lemma~\ref{AAAA} the following cases for $(\fh_j,V)$ need to be checked:
\begin{enumerate}
\item[-] $(\fh_j,\omega_1)$ where $\fh_j$ is one of $\sl(n)$, $\sp(n)$, $\so(n)$,
  $\ssG_2$, or $\ssE_6$;
\item[-] $(\sl(n),\omega_2)$;
\item[-] the spin representation for $n=7,8,9,10$.
\end{enumerate}

In the cases marked ``$b$'', the dimension of the Borel subalgebra of
$\tilde\fl$ is smaller than $\|dim| V$. The same happens for ``$b'$'' but
one has to take the image of $\fl$ in $\fh_j$ into account. In the
cases marked ``$+$'', $V$ contains two $\fl$\_stable lines with the
same character, hence cannot be $\tilde\fl$\_spherical. The
cases marked by an ``$x$'' will be checked separately:

$(\ssE_6, \ssF_4)$ is not in Table~\ref{PairsForTriples} because $\ssF_4$
does not have non-trivial spherical modules.

For $(\ssF_4,\ssB_4)$, the principal subalgebra is $\fl = \ssB_3$,
spin\_embedded into $\mf{so}(8) \inn \ssB_4$. $(\ssB_4, \omega_4)$ is
not spherical for $\tilde\fl$ because its Borel subalgebra is too
small; and $(\ssB_4,\omega_1)$ is not because the orbits of
$\mf{spin}(7)$ have codimension at least $2$.

For $(\ssE_6,\ssD_5)$, the principal subalgebra is $\fl =
\ssA_3+\ft^1$ where $\ssA_3 \into \ssB_3 =\mf{spin}(7) \into \ssD_4
\inn \ssD_5$. Therefore, the same reasoning as in the case $(\ssF_4,
\ssB_4)$ applies.

The last $x$-case is $(\sl(m{+}2)+\sp(2n{+}2),\sl(m)+\sl(2)+\sp(2n))$.
The principal subalgebra is $\fl=\gl(m{-}2)+\sp(2n{-}2)$. We see that
$(\gl(m),\omega_1)$ and $(\sp(2n),\omega_1)$ are not spherical for
$\tilde\fl$ by ``$+$'', and that $(\gl(m),\omega_2)$ is not by ``$b$''.

Finally, we come the two cases which are marked ``$\circ$''. These
form series which consist partially of possible base components for
triples with nonzero fiber. More precisely,

\begin{enumerate}
\item Let $m\ge n\ge 1$ and $m+n\ge 7$. If the spherical pair $(\so(m+n),
  \so(m)+\so(n))$ is a base component then $n=1$ or $n=2$.\label{soso}
\item Let $\fg$ be a simple Lie algebra. If the spherical pair
  $(\fg+\fg, \fg)$ is a base component then $\fg \isom \sl(n)$, $n
  \ge2$.\label{gggg}
\end{enumerate}

\medskip
\noindent{\it Proof of \REF{soso}.}
Suppose $n\ge3$. Then $\fl=\so(m-n)\subseteq\so(m)$ implies that $V$ is
a spherical $\so(m)$\_module. Inspecting Table~\ref{sphermod} we see that we
have to check only for $V$ the defining or the spin representation
(for $m=4$ and $m=6$ we have to use Lemma~\ref{AAAA}\REF{pleth}).
Both acquire multiplicities when restricted to $\fl=\so(m-n)$.

\medskip\noindent{\it Proof of \REF{gggg}.}
  The principal subalgebra is a Cartan subalgebra of $\fg$. Thus, in
  order for $(\fg+\fg,\fg,V)$ to be spherical we must have $\dim
  V\le\mathop{\rm rk}\fg+1$. Inspecting Table~\ref{sphermod}, we see this is
  only possible for $\fg=\sl(n)$.
\medskip

This finishes the proof of Proposition~\ref{ActualComp}.
\end{proof}

\section{Reduced triples}

In this section, we justify the inference rules.

\begin{lemma}
\label{nonred1}

  Fix $m\ge1$. Let $\fh_0\subseteq\fg_0$ be semisimple Lie algebras
  and let $V$ be an $\sl(2)+\fh_0$\_module. Then the following
  are equivalent:

\begin{enumerate}
\item $(\tilde\fg,\tilde\fh,V)=(\sl(2)+\fg_0,\sl(2)+\fh_0,V)$ is
a primitive spherical triple.

\item $(\fg,\fh,V)=(\sp(2m+2)+\fg_0,\sp(2m)+\sl(2)+\fh_0,V)$ is a
primitive spherical triple. 
\end{enumerate}
\end{lemma}

\begin{proof}
The principal subalgebra of $(\sp(2m+2),\sp(2m)+\sl(2))$ surjects onto
$\sl(2)$. The statement now follows from Theorem~\ref{crit}.
\end{proof}

\begin{lemma} Let $\fh_0\subseteq\fg_0$ be semisimple Lie
  algebras and let $V$ be an $\so(6)+\fh_0$\_module. The following
  are equivalent:\label{so5so6}

\begin{enumerate}
\item\label{so6} $(\fg,\fh,V)=(\so(7)+\fg_0,\so(6)+\fh_0,V)$ is
a primitive spherical triple.

\item\label{so5} $(\tilde\fg,\tilde\fh,V)=(\so(5)+\fg_0,\so(5)+\fh_0,V)$ is
a primitive spherical triple such that $V|_{\so(5)}$ contains only the
trivial and the spin represntation.

\end{enumerate}

\noindent
Moreover, $V$ contains in this case, as an $\so(6)$-module, at most the
trivial and the two spin representations. All of these stay
irreducible as $\so(5)$\_modules.
\end{lemma}

\begin{proof}
  First assume \REF{so6}. The principal subalgebra of
  $(\so(7),\so(6))$ is $\so(5)$. Thus, Theorem~\ref{crit} implies that
  $(\tilde\fg,\tilde\fh,V)=(\so(5)+\fg_0,\so(5)+\fh_0,V)$ is
  spherical. Let $U$ be a simple component of $V|_{\so(6)}$. Then
  the erasing lemma (Corollary~\ref{erase}) implies that
  $(\so(7),\so(6),U)$ is spherical. Table~\ref{sphermod} implies easily
  that $U$ is either $\CC$, $\CC^4$, or $(\CC^4)^*$. This proves
  \REF{so5} and the last statement.
  
  Now assume \REF{so5}. It is not possible to deduce the sphericality
  of $(\fg,\fh,V)$ directly from that of $(\tilde\fg,\tilde\fh,V)$. A
  counterexample is $(\so(7),\so(6),\CC^6)$. The problem is that the
  torus $\fz$ may be different. On the other hand, the assumption on
  $V|_{\so(5)}$ implies easily that $V|_{\so(6)}$ contains only the
  trivial and the spin representations. This means that $V|_{\so(5)}$
  and $V|_{\so(6)}$ have the same number of irreducible components,
  i.e., $\fz$ does not change.
\end{proof}

\begin{defn}
A spherical triple $(\fg,\fh,V)$ is called {\it reduced} if it is
primitive and
\begin{enumerate}
\item $(\so(7),\so(6))$ is not a component of $(\fg,\fh)$ and
\item $(\sp(2m+2),\sp(2m)+\sl(2))$ is, for any $m\ge1$, not a
  component of $(\fg,\fh)$ such that $\sp(2m)$ acts trivially on $V$.
\end{enumerate}
\end{defn}

\begin{cor}
All primitive spherical triples can be obtained from reduced spherical
triples by (possibly repeated) application of the inference rules.
\end{cor}

Thus, Tables~\ref{sphermod} and~\ref{newtab2} constitute, in fact, a
classification of reduced spherical triples.

\section{Simple extensions}

As one sees from the tables, the bulk of the new spherical triples is
of the following type.

\begin{defn}
A reduced spherical triple $(\fg,\fh,V)$ is called a {\it simple
  extension (of $(\fg_0,\fh_0)$)} if

\begin{enumerate} 
\item it has exactly one base component $(\fg_0,\fh_0)$,
\item it has exactly one fiber component, and
\item the intersection of the base component with the fiber component 
is a single vertex $\fh_j$.
\end{enumerate}
\end{defn}

In other words, the diagram of $(\fg,\fh,V)$ is obtained by gluing a
diagram of Table~\ref{PairsForTriples} to a diagram of
Table~\ref{sphermod} to one vertex in the $\fh$\_level. In this
section we prove:

\begin{prop} All simple extensions are contained in
  Table~\ref{newtab2}.
\end{prop}

We start by determining all simple extensions which are glued at an $\sl(2)$:

\begin{lemma}
\label{sl2glue}
Let $(\fg,\fh,V)$ be a simple extension glued at
$\fh_j=\sl(2)$. Then $(\fg,\fh,V)$ is obtained by gluing any of the
components
\begin{enumerate}
\item[---] $(\sl(n+2),\underline{\sl(2)}+\sl(n)), n\ge1$
\item[---] $(\sp(4)+\sp(2n+2),\underline{\sl(2)}+\sl(2)+\sp(2n)), n\ge0$
\item[---] $(\sp(2m+2)+\sp(2n+2),\sp(2m)+\underline{\sl(2)}+\sp(2n)), m,n\ge0$
\end{enumerate}
with any of the modules
\begin{enumerate}
\item[---] $(\underline{\sl(2)}+\sl(n), \omega_1\omega_1')$, $n \ge 1$
\item[---] $(\underline{\sl(2)}+\sp(2n), \omega_1\omega_1')$, $n \ge 2$
\end{enumerate}
at an underlined factor. All of these are contained in Table~\ref{newtab2}.
\end{lemma}

\begin{proof}
Let $\fl_0$ be the image of $\fl\subseteq\fh$ in $\fh_j$. Then there
are three possibilities: (1) $\fl_0=0$. In this case, one cannot glue
$\fh_j$ to anything because, as an $\fl_0+\fz$\_module, $V$ would have
multiplicities. (2) $\fl_0=\sl_2$. Inspecting Table~\ref{PairsForTriples}
this is precisely the case for
$(\sp(2m+2),\sl(2)+\sp(2m))$ which we dealt with in
Lemma~\ref{nonred1}. (3) $\fl_0=\ft^1$. These are precisely the cases
listed in the lemma.

Now we claim that the listed modules are those which stay spherical
when $\sl(2)$ is replaced by $\ft^1$. In fact, a spherical module
$(\sl(2)+\fh_0,V)$ has this property if and only of the triple
$(\sl(3)+\fh_0,\sl_2+\fh_0,V)$ is spherical. Now we can apply
Lemma~\ref{reductionlemma} with $U=\omega_1$. Thus, we see that the
triple is spherical if and only if
$(\sl(2)+\fh_0,\sl(2)+\fh_0,U\oplus V)$ is spherical. Using
Table~\ref{sphermod} we easily obtain the list of modules given.
\end{proof}

For the rest of the proof we are going through the list of possible
base components (Table~\ref{PairsForTriples}). The underlined
factor is the one we want to glue, in case there is a choice.

\subsection{$(\fg_0,\fh_0)=(\sl(m+n),\underline{\sl(m)}+\sl(n))$,
  $m\ge3$, $n\ge1$:}\label{trick}
In this case Lemma~\ref{reductionlemma} applies with
$U=\omega_1\omega_1'$. Moreover, the action of $\bar\fh$ on $U$
contains the scalars. This implies that we may replace
$(\sl(m+n)+\tilde\fg,\sl(m)+\sl(n)+\tilde\fh,V)$ by
$(\sl(m)+\sl(n)+\tilde\fg,\sl(m)+\sl(n)+\tilde\fh,\omega_1\omega_1'\oplus
V)$. Graphically:
$$
\begin{texdraw}
\labelL(0,1)[0]{\sl(m)}
\strichRO(0,1)
\strichU(0,1)
\labelO(1,2)[0]{\sl(m+n)}
\strichRU(1,2)
\labelR(2,1)[0]{\sl(n)}
\strichRU(2,1)
\end{texdraw}
\qquad\raise20pt\hbox{$\Longleftrightarrow$}\qquad
\begin{texdraw}
\labelL(0,1)[0]{\sl(m)}
\strichRU(0,1)
\strichU(0,1)
\punkt(1,0)
\strichRO(1,0)
\labelR(2,1)[0]{\sl(n)}
\strichRU(2,1)
\end{texdraw}
$$
Taking also the degenerate case $n=1$ into account we see from
Table~\ref{sphermod} that only the triples
$(\sl(m+1)+\sl(k),\sl(m)+\sl(k),\omega_1\omega_1')$, $k\ge1$,
$(\sl(m+1),\sl(m),\omega_2)$, and $(\sl(m+n),\sl(m)+\sl(n),\omega_1)$
are spherical (for $m\ge3$).

\subsection{$(\fg_0,\fh_0)=(\sl(2n),\sp(2n))$, $n\ge2$:}\label{slsp}
we have $\fl=\sl(2)^n$ and we can glue the following modules:

\subsubsection{$(\sp(2n),\omega_1)$:} the glued triple is indeed spherical.

\subsubsection{$(\sp(4),\omega_2)$ with $n=4$:} the dimension of $V$
is too big.

\subsubsection{$(\sp(2n)+\sl(m),\omega_1\omega_1')$ with $n\ge3$,
  $m\ge2$:} the module $(\fl,V)$ is not in Table~\ref{sphermod}.

\subsubsection{$(\sp(4)+\sl(m),\omega_1\omega_1')$ with $n=4$,
  $m\ge2$:} if $m\ge3$, $(\fl,V)$ is not in Table~\ref{sphermod}. For
$m=2$, the dimension of $V$ is too big.
\medskip

This exhausts all irreducible $V$. The only case left to check is

\subsubsection{$(\sp(2n),\omega_1+\omega_1)$:} the dimension of $V$
is too big.

\subsection{$(\fg_0,\fh_0)=(\sp(2m+2n),\underline{\sp(2m)}+\sp(2n))$, $m\ge2,
n\ge1$:}\label{spspsp} if $m\le n+1$ then the image of $\fl$ in
$\sp(2m)$ is $\sl(2)^m$ and we are back to case~\ref{slsp}
with only $(\sp(2m+2n),\sp(2m)+\sp(2n),\omega_1)$ being
spherical. Therefore assume $m>n+1$. Since, in particular, $m\ge3$ we
have to check the following cases:

\subsubsection{$(\sp(2m),\omega_1)$:} spherical.

\subsubsection{$(\sp(2m)+\sl(k),\omega_1\omega_1')$ with $k=2,3$:}
here $(\fl,V)$ is not in Table~\ref{sphermod} unless $n=1$ and
$k=2$. Also in that case, $V$ is not a spherical $\fl+\fz$\_module (see
Lemma~\ref{deficient}).

\medskip
After this, the only reducible $V$ to check is

\subsubsection{$(\sp(2m),\omega_1+\omega_1)$:} not spherical (see
Lemma~\ref{deficient}).

\subsection{$(\fg_0,\fh_0)=(\so(n+1),\so(n))$, $n\ge7$:} in this case
$\fl=\so(n-1)$ and the following possible gluings arise:

\subsubsection{$(\so(n),\omega_1)$:} $V|_{\fl+\fz}=\CC^{n-1}\oplus\CC$
is not spherical.

\subsubsection{$(\so(n),\hbox{\rm spin rep.})$:} according to
Table~\ref{sphermod} we must have $n=7$, $8$, $9$, or $10$. All of
them yield spherical triples except for $n=9$ (see
Lemma~\ref{deficient}).
\medskip
As for reducible representations we have:

\subsubsection{$(\so(8),\omega_3+\omega_4)$:} not spherical since
$V|_\fl$ has multiplicities.

\subsection{$(\fg_0,\fh_0)=(\so(n+2),\so(n))$, $n\ge5$:} We can use
Lemma~\ref{reductionlemma} with $U=\omega_1$. This implies $n=8$ and
$V=\hbox{spin rep.}$

\subsection{$(\fg_0,\fh_0)=(\so(2n),\sl(n))$, $n\ge4$:} using
Lemma~\ref{reductionlemma} with $U=\omega_2$ we see that only
$(\so(2n),\sl(n),\omega_1)$ is spherical.

\subsection{$(\fg_0,\fh_0)=(\ssG_2,\sl(3))$:} in this case we can apply
Lemma~\ref{reductionlemma} with $U=\omega_1$. Since $\fs=0$ only the
implication $(\ref{reditem1})\Rightarrow(\ref{reditem2})$ is valid. This leaves the following
cases to check:

\subsubsection{$(\sl(3),\omega_1)$:} spherical.

\subsubsection{$(\sl(3)+\sl(n),\omega_1\omega_1')$, $n\ge2$:} not
spherical (see Lemma~\ref{deficient}).

\subsection{$(\fg_0,\fh_0)=(\sl(n)+\sl(n),\sl(n))$, $n\ge3$:} when only
$\sl(n)$ acts on $V$ then $\fl=\ft^{n-1}$, hence $\dim V\le n$. This
implies $(\fh,V)=(\sl(n),\omega_1)$ which indeed is spherical for
$\ft^{n-1}+\ft^1$. Otherwise, $V$ has an irreducible component of the
form $(\sl(n)+\fk,\omega_1\otimes U)$ with $\|dim|U\ge2$ such that
$\CC^n\otimes U$ is spherical for $\ft^n+\fk$. This does not exist
for $n\ge3$ as Table~\ref{sphermod} shows there are no indecomposable spherical modules with more than two irreducible (fiber) components.

\subsection{$(\fg_0,\fh_0)=(\sp(2m+2)+\sp(2n+2),
  \underline{\sp(2m)}+\sl(2)+\sp(2n))$, $m\ge2$, $n\ge0$:} after
cutting the diagram at the $\sl(2)$\_ vertex of $\fh$
(Lemma~\ref{cutting}) we see from case~\ref{spspsp} that at most the
triple $(\sp(2m+2)+\sp(2n+2),
\sp(2m)+\sl(2)+\sp(2n),\omega_1)$ is spherical, and it is.

\section{Tree-like extensions}

In this section, we classify all reduced spherical triples whose
diagram is a tree, i.e., contains no cycles. These triples are called
{\it tree\_like}.

\begin{lemma}
\label{onecomp}
  Let $(\fg,\fh,V)$ be a reduced tree\_like spherical triple with
  exactly one base component. Then it appears in
  Table~\ref{newtab2}.
\end{lemma}

\begin{proof}
Let $(\fg_0,\fh_0)$ be the base component.  We handle the case of
exactly two fiber components first. Since they have to be glued to
different simple factors of $\fh_0$ we see that $\fh_0$ cannot be
simple. Moreover, $(\fg_0,\fh_0)$ must have at least one simple
extension. That leaves the following possibilities:

\subsection{$(\fg_0,\fh_0)=(\sl(m+n),\sl(m)+\sl(n))$, $m,n\ge2$:}
Again, we can use Lemma~\ref{reductionlemma}. Then we are done, since
there are no indecomposable spherical modules with more than two
irreducible components.

\subsection{$(\fg_0,\fh_0)=(\sp(2m+2n),\sp(2m)+\sp(2n))$,
  $m,n\ge1$:} Let $V=V'\oplus V''$ where $V'$ and
  $V''$ are attached to $\sp(2m)$ and $\sp(2n)$, respectively.

\subsubsection{Case $m\ge n\ge2$:} then $V'=\omega_1$ and
$V''=\omega_1'$. Hence, $V$ contains as an $\fl$\_submodule
$(\sl(2)^n,(\CC^2+\CC^2)^n)$ which is not spherical for dimension
reasons.

\subsubsection{Case $m>n=1$:} then
$V'=\omega_1=\CC^{2m-2}\oplus\CC^2$. The $\sl(2)$\_factor in $\fl$
acts diagonally on $\CC^2\subseteq V'$ and on $V''$. Thus, we can
attach a representation $V''$ if and only if $\CC^2\oplus V''$ is a
spherical $\fl+\fz$\_module. These are precisely the representations
$(\sl(2)+\sl(n),\omega_1\omega_1')$, $n\ge1$ and
$(\sl(2)+\sp(2n),\omega_1\omega_1')$, $n\ge2$ (as in the proof of
Lemma~\ref{sl2glue}).

\subsubsection{Case $m=n=1$:} in this case $\fl=\sl(2)$ acting
diagonally on $V'$ and $V''$. Therefore, $V'$ and $V''$ are two
indecomposable modules which, when branched to $\fl$, are ``glued'' at one $\sl(2)$\_vertex
yielding the last item of Table~\ref{newtab2}.

\subsection{$(\fg_0,\fh_0)=(\sp(2m+2)+\sp(2n+2),\sp(2m)+\sl(2)+\sp(2n))$,
  $m,n\ge0$:} the principal subalgebra of the pair contains a factor
$\ft^1$. Looking at all simple extensions one sees that each
irreducible component of $V$ contains an $\fl$\_submodule of the form
$\CC^2\otimes\CC^k$, $k\ge1$, where $\ft^1\subseteq\fl$ acts only via
its embedding into $\sl(2)$ and where $\CC^k$ is acted on by either $\sl(k)$
or $\sp(k)$. Since there are two fiber components we have another
submodule $\CC^2\otimes\CC^l$, $l\ge1$. The $\ft^1$\_factor acts
diagonally on the two submodules of $V$. Since
$(\CC^2\otimes\CC^k)\oplus(\CC^2\otimes\CC^l)$ is not a spherical
module for $\ft^1+\sl(k)+\sl(l)+\fz$ we conclude that $(\fg_0,\fh_0)$
has no extensions with two fiber components.

\medskip
Finally, none of the triples with two fiber components has a ``free''
factor in $\fh_0$. Therefore, there are no reduced triples with more
than two fiber components.
\end{proof}

\begin{lemma}
\label{twocomp}
  There are no reduced tree\_like spherical triples with two base
  components.
\end{lemma}

\begin{proof}
  Let $(\fg,\fh,V)$ be a counterexample whose graph $\Gamma$ has a
  minimal number of edges. Since $\Gamma$ is connected there must be a
  fiber\_component $F$ which is connected to two base components
  $B_1$ and $B_2$. One may erase the edges leading
  to other components (Corollary~\ref{erase}). Hence minimality implies
  that $\Gamma$ has only the three components $B_1, B_2$ and $F$. Since $\Gamma$ is a
  tree the intersection $F\cap B_i$ consists of a single vertex $\fh_i$.
  
  Let $\Gamma_i$ be the union of $F$ with $B_i$. Then cutting $\Gamma$ at
  $\fh_{3-i}$ (Lemma~\ref{cutting}) shows that $\Gamma_i$ is a
  connected component of a spherical triple hence spherical as well.
  
  Lemma~\ref{onecomp} says that we can find $\Gamma_i$ in
  Table~\ref{newtab2}. Since exactly two simple factors, $\fh_1$ and
  $\fh_2$, act on $V$ we see that $V$ is either $(\sl(2)+\sp(2n),\omega_1\omega_1')$
  with $n\ge2$ or $(\sl(m)+\sl(n),\omega_1\omega_1')$ with $m,n\ge2$. We can rule
  out the first case since the $\sp(2n)$\_factor cannot be attached to
  anything. When $m,n\ge3$, the second case leads to a unique
  triple namely
$$
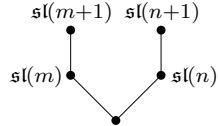
\begin{texdraw}
\strichU(0,2)
\strichU(2,2)
\strichRU(0,1)
\strichRO(1,0)
\punkt(1,0)
\labelO(0,2)[0]{\sl(m{+}1)}
\labelO(2,2)[0]{\sl(n{+}1)}
\labelL(0,1)[0]{\sl(m)}
\labelR(2,1)[0]{\sl(n)}
\end{texdraw}
$$
For $m,n\ge2$ none of these triples is spherical. If $m=2$ or
$n=2$ other base components (namely those from Lemma~\ref{sl2glue}) could
be attached. All of them just act through a factor $\ft^1$ in
$\fl$. Therefore, these triples are not spherical either.
\end{proof}

Combining Lemmas~\ref{onecomp} and~\ref{twocomp} we get:

\begin{cor} Let $(\fg,\fh,V)$ be a reduced tree\_like spherical
  triple with $\fg\ne\fh$ and $V\ne0$. Then it appears in
  Table~\ref{newtab2}.
\end{cor}

\section{Cycles}

The following lemma finishes off the classification:

\begin{lemma}
Every spherical triple is tree\_like.
\end{lemma}

\begin{proof}
  Let $(\fg,\fh,V)$ be a counterexample whose diagram $\Gamma$ has
  minimal number of edges. Clearly, the triple is reduced. Let $C$ be
  a cycle. Since, by inspection, all base components
  $(\fg_0,\fh_0)$ are trees, $C$ has to contain a vertex $V_k$.
  Moreover, again by inspection, at most, hence exactly two factors of
  $\fh$ act non\_trivially on $V_k$. Thus $V_k=U_1\otimes U_2$. Let
  $\tilde\Gamma$ be the triple obtained by cutting $\Gamma$ at $V_k$
  (Lemma~\ref{Vtensor}), i.e., where we replace $U_1\otimes U_2$ by
  $U_1\oplus U_2$. Then $\tilde\Gamma$ is a tree since otherwise we
  could erase the edges adjacent to $U_i$ and obtain a smaller spherical
  triple whose diagram contains a cycle. This implies that
  $\tilde\Gamma$ is a member of Table~\ref{newtab2}. Moreover, it
  contains two vertices in the $V$\_layer which are adjacent to exactly
  one vertex each in the $\fh$\_layer. There is only one such case which
  results in
$$
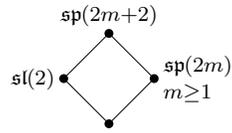
\begin{texdraw}
\strichRO(0,1)
\strichRU(1,2)
\strichRO(1,0)
\strichRU(0,1)
\labelL(0,1)[0]{\sl(2)}
\labelO(1,2)[0]{\sp(2m+2)}
\labelR(2,1)[.3]{\sp(2m)}
\labelRR(2,.7)[0]{m\ge1}
\punkt(1,0)
\end{texdraw}
$$
for $\Gamma$. As a module for $\tilde\fl=\sp(2)+\sp(2m-2)+\ft^1$, $V$
contains $(\gl(2),\CC^2\otimes\CC^2)$ which is not spherical, by
dimension.
\end{proof}

\newpage

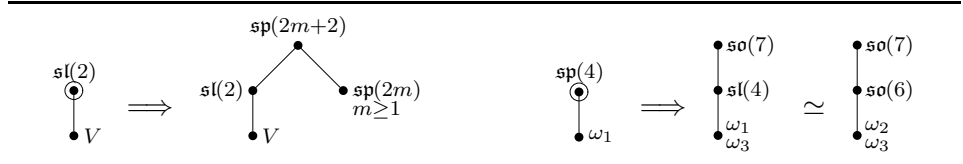
\begin{table}[htbp]
\caption{Primitive spherical modules}
\label{sphermod}
\hrule
\vskip5pt
\baselineskip30pt
\leavevmode
\begin{texdraw}
\strichRU(0,1)
\strichRO(1,0)
\labelO(0,1)[0]{m\ge1}
\labelOO(0,1.4)[0]{\sl(m)}
\labelO(2,1)[0]{n\ge2}
\labelOO(2,1.4)[0]{\sl(n)}
\punkt(1,0)
\kreis(0,1)
\kreis(2,1)
\end{texdraw}
%
%
%
\hskip0pt plus 1fil
\begin{texdraw}
\strichU(0,1)
\labelR(0,1)[.2]{\sl(n)}
\labelRR(0,.6)[.2]{n\ge2}
\labelR(0,0)[0]{\omega_1^2}
\kreis(0,1)
\end{texdraw}
\hskip0pt plus 1fil
\begin{texdraw}
\strichU(0,1)
\labelR(0,1)[.2]{\sl(n)}
\labelRR(0,.6)[.2]{n\ge4}
\labelR(0,0)[0]{\omega_2}
\end{texdraw}
\hskip0pt plus 1fil
\begin{texdraw}
\strichU(0,1)
\labelR(0,1)[.2]{\sp(2n)}
\labelRR(0,.6)[.2]{n\ge2}
\punkt(0,0)
\kreis(0,1)
\end{texdraw}
\hskip0pt plus 1fil
\begin{texdraw}
\strichU(0,1)
\labelR(0,1)[0]{\sp(4)}
\labelR(0,0)[0]{\omega_2}
\end{texdraw}
\hskip0pt plus 1fil
\begin{texdraw}
\strichRU(0,1)
\strichRO(1,0)
\labelO(0,1)[0]{n\ge2}
\labelOO(0,1.4)[0]{\sp(2n)}
\labelO(2,1)[0]{\sl(2)}
\punkt(1,0)
\kreis(0,1)
\kreis(2,1)
\end{texdraw}
\hskip0pt plus 1fil
\begin{texdraw}
\strichRU(0,1)
\strichRO(1,0)
\labelOO(0,1.4)[0]{\sp(2n)}
\labelO(0,1)[0]{n\ge3}
\labelO(2,1)[0]{\sl(3)}
\punkt(1,0)

\end{texdraw}
\hskip0pt plus 1fil
\begin{texdraw}
\strichRU(0,1)
\strichRO(1,0)
\labelO(0,1)[0]{\sp(4)}
\labelO(2,1)[0]{n\ge3}
\labelOO(2,1.4)[0]{\sl(n)}
\punkt(1,0)
\kreis(0,1)
\end{texdraw}
\hskip0pt plus 1fil
\begin{texdraw}
\strichU(0,1)
\labelR(0,1)[.2]{\so(n)}
\labelRR(0,.6)[.2]{n\ge7}
\punkt(0,0)
\end{texdraw}
\hskip0pt plus 1fil
\begin{texdraw}
\strichU(0,1)
\labelR(0,1)[0]{\so(7)}
\labelR(0,0)[0]{\omega_3}
\end{texdraw}
\hskip0pt plus 1fil
\begin{texdraw}
\strichU(0,1)
\labelR(0,1)[0]{\so(9)}
\labelR(0,0)[0]{\omega_4}
\end{texdraw}
\hskip0pt plus 1fil
\begin{texdraw}
\strichU(0,1)
\labelR(0,1)[0]{\so(10)}
\labelR(0,0)[0]{\omega_5}
\end{texdraw}
\hskip0pt plus 1fil
\begin{texdraw}
\strichU(0,1)
\labelR(0,1)[0]{\ssG_2}
\labelR(0,0)[0]{\omega_1}
\end{texdraw}
\hskip0pt plus 1fil
\begin{texdraw}
\strichU(0,1)
\labelR(0,1)[0]{\ssE_6}
\labelR(0,0)[0]{\omega_1}
\end{texdraw}
\hskip0pt plus 1fill\

\vskip10pt
\hrule
\vskip5pt

\leavevmode\noindent
\begin{texdraw}
\strichRU(0,1)
\strichRO(1,0)
\strichRU(2,1)
\labelO(0,1)[0]{m\ge1}
\labelOO(0,1.4)[0]{\sl(m)}
\kreis(0,1)
\punkt(1,0)
\labelO(2,1)[0]{n\ge3}
\labelOO(2,1.4)[0]{\sl(n)}
\labelRR(3,0)[.2]{\omega_1}
\labelR(3,0)[-.2]{\omega_{n{-}1}}
\end{texdraw}
\hskip0pt plus 1fil
\begin{texdraw}
\strichRO(0,0)
\strichRU(1,1)
\labelO(1,1)[0]{n\ge4}
\labelOO(1,1.4)[0]{\sl(n)}
\labelR(0,0)[0]{\omega_2}
\labelRR(2,0)[.2]{\omega_1}
\labelR(2,0)[-.2]{\omega_{n{-}1}}
\end{texdraw}
\hskip0pt plus 1fil
\begin{texdraw}
\strichRO(0,0)
\strichRU(1,1)
\labelO(1,1)[0]{n\ge2}
\labelOO(1,1.4)[0]{\sp(2n)}
\kreis(1,1)
\punkt(0,0)
\punkt(2,0)
\end{texdraw}
\hskip0pt plus 1fil
\begin{texdraw}
\strichRO(0,0)
\strichRU(1,1)
\labelO(1,1)[0]{\so(8)}
\labelR(0,0)[0]{\omega_3}
\labelR(2,0)[0]{\omega_4}
\end{texdraw}
\hskip0pt plus 1fil
%
\begin{texdraw}
\strichRU(0,1)
\strichRO(1,0)
\strichRU(2,1)
\strichRO(3,0)
\labelLL(0,1)[.4]{\sl(m)}
\labelL(0,1)[0]{\sp(2m)}
\labelLL(0,1)[-.4]{m\ge1}
\kreis(0,1)
\labelO(2,1)[0]{\sl(2)}
\kreis(2,1)
\labelRR(4,1)[.4]{\sl(n)}
\labelR(4,1)[0]{\sp(2n)}
\labelRR(4,1)[-.4]{n\ge1}
\kreis(4,1)
\punkt(1,0)
\punkt(3,0)
\end{texdraw}
\hskip0pt plus 1fill\

\vskip5pt

\hrule
\end{table}

\begin{table}[htbp]
\caption{Reduced spherical triples $(\fg,\fh,V)$, $\fg\ne\fh$, $V\ne0$}
\label{newtab2}
\hrule
\vskip5pt
\baselineskip50pt
\leavevmode
\begin{texdraw}
\strichRO(0,1)
\strichRU(1,2)
\strichU(2,1)
\punkt(2,0)
\labelUu(0,1)[0]{\sl(m)}
\labelUU(0,.6)[0]{m \ge 2}
\labelO(1,2)[0]{\sl(m{+}n)}
\labelR(2,1)[.3]{\sl(n)}
\labelRR(2,.5)[.3]{n\ge 3}
\end{texdraw}
\hskip0pt plus 1fil
\begin{texdraw}
\strichRO(0,1)
\strichRU(1,2)
\strichRU(2,1)
\strichRO(3,0)
\punkt(3,0)
\kreis(4,1)
\labelUu(0,1)[0]{\sl(m)}
\labelUU(0,.6)[0]{m \ge 1}
\labelO(1,2)[0]{\sl(m{+}2)}
\labelU(2,1)[-.4]{\sl(2)}
\labelR(4,1)[.5]{\sl(n)}
\labelRR(4,.5)[.5]{\sp(2n)}
\labelRR(4,0)[.5]{n \ge 1}
\end{texdraw}
\hskip0pt plus 1fil
\begin{texdraw}
\strichU(0,2)
\strichRU(0,1)
\strichRO(1,0)
\punkt(1,0)
\kreis(2,1)
\labelO(0,2)[0]{\sl(m{+}1)}
\labelUu(0,1)[-.6]{\sl(m)}
\labelUU(0,.6)[-.6]{m \ge 3}
\labelR(2,1)[.3]{\sl(n)}
\labelRR(2,.5)[.3]{n\ge 1}
\end{texdraw}
\hskip0pt plus 1fil
\begin{texdraw}
\strichU(0,2)
\strichU(0,1)
\labelR(0,2)[0]{\sl(n{+}1)}
\labelR(0,1)[0.3]{\sl(n)}
\labelRR(0,0.6)[.3]{n \ge 4}
\labelR(0,0)[0]{\omega_2}
\end{texdraw}
\hskip0pt plus 1fil
\begin{texdraw}
\strichRU(0,2)
\strichRO(1,1)
\strichU(1,1)
\punkt(1,0)
\labelO(0,2)[0]{\sl(n)}
\labelO(2,2)[0]{\sl(n)}
\labelR(1,1)[-.3]{\sl(n)}
\labelRR(1,.5)[-.3]{n \ge 2}
\end{texdraw}
\hskip0pt plus 1fil
\begin{texdraw}
\strichU(0,2)
\strichU(0,1)
\punkt(0,0)
\labelR(0,2)[0]{\sl(2n)}
\labelR(0,1)[0.3]{\sp(2n)}
\labelRR(0,.5)[0.3]{n \ge 2}
\end{texdraw}
\hskip0pt plus 1fil
\begin{texdraw}
\strichU(0,2)
\strichU(0,1)
\punkt(0,0)
\labelR(0,2)[0]{\so(2n)}
\labelR(0,1)[0.3]{\sl(n)}
\labelRR(0,.5)[0.3]{n \ge 4}
\end{texdraw}
\hskip0pt plus 1fil
\begin{texdraw}
\strichU(0,2)
\strichU(0,1)
\labelR(0,2)[0]{\so(n{+}1)}
\labelR(0,1)[0.3]{\so(n)}
\labelRR(0,0.6)[0.3]{n=7,8,10}
\labelR(0,0)[0]{{\rm spinrep.}}
\end{texdraw}
\hskip0pt plus 1fil
\begin{texdraw}
\strichU(0,2)
\strichU(0,1)
\labelR(0,2)[0]{\so(10)} 
\labelR(0,1)[0]{\so(8)} 
\labelR(0,0)[0]{{\rm spinrep.}}
\end{texdraw}
\hskip0pt plus 1fil
\begin{texdraw}
\strichU(0,2)
\strichU(0,1)
\punkt(0,0)
\labelR(0,2)[0]{\ssG_2}
\labelR(0,1)[0]{\sl(3)}
\end{texdraw}
\hskip0pt plus 1fil
\begin{texdraw}
\strichRO(0,1)
\strichRU(1,2)
\strichU(2,1)
\punkt(2,0)
\labelUu(0,1)[0]{\sp(2m)}
\labelUU(0,.5)[0]{m\ge1}
\labelO(1,2)[0]{\sp(2m{+}2n)}
\labelR(2,1)[0.3]{\sp(2n)}
\labelRR(2,.5)[0.3]{n\ge1}
\end{texdraw}
\hskip0pt plus 1fil
\begin{texdraw}
\strichRO(0,1)
\strichRU(1,2)
\strichRO(2,1)
\strichRU(3,2)
\strichU(4,1)
\punkt(4,0)
\labelUu(0,1)[0]{\sp(2m)}
\labelUU(0,.5)[0]{m\ge0}
\labelO(1,2)[-.1]{\sp(2m{+}2)}
\labelU(2,1)[0]{\sl(2)}
\labelO(3,2)[0.1]{\sp(2n{+}2)}
\labelR(4,1)[0.3]{\sp(2n)}
\labelRR(4,.5)[0.3]{n\ge2}
\end{texdraw}
\hskip0pt plus 1fil
\begin{texdraw}
\strichRO(0,1)
\strichRU(1,2)
\strichRO(2,1)
\strichRU(3,2)
\move(2 1)\rlvec(2 -1)
\move(4 0)\rlvec(2 1)
\punkt(4,0)
\kreis(6,1)
\labelUu(0,1)[0]{\sp(2l)}
\labelUU(0,.5)[0]{l\ge0}
\labelO(1,2)[0]{\sp(2l{+}2)}
\labelU(2,1)[-.4]{\sl(2)}
\labelO(3,2)[.2]{\sp(2m{+}2)}
\labelR(4,1)[.4]{\kern-5pt\sp(2m)}
\labelRR(4,.5)[.4]{m\ge0}
\labelR(6,1)[.5]{\sl(n)}
\labelRR(6,.5)[.5]{\sp(2n)}
\labelRR(6,0)[.5]{n\ge1}
\end{texdraw}
\hskip0pt plus 1fil
\begin{texdraw}
\strichRO(0,1)
\strichRU(1,2)
\strichRO(2,1)
\strichRU(3,2)
\strichRU(4,1)
\strichRO(5,0)
\punkt(5,0)
\kreis(6,1)
\labelUu(0,1)[0]{\sp(2m)}
\labelUU(0,.5)[0]{m\ge0}
\labelO(1,2)[0]{\sp(2m{+}2)}
\labelU(2,1)[0]{\sl(2)}
\labelO(3,2)[0]{\sp(4)}
\labelU(4,1)[-.4]{\sl(2)}
\labelRR(6,1.5)[0]{\sl(n)}
\labelR(6,1)[0]{\sp(2n)}
\labelRR(6,.5)[0]{n\ge1}
\end{texdraw}
\hskip0pt plus 1fil
\begin{texdraw}
\strichU(0,1)
\strichRO(0,1)
\strichRU(1,2)
\strichRU(2,1)
\strichRO(3,0)
\punkt(0,0)
\punkt(3,0)
\kreis(4,1)
\labelL(0,1)[.3]{\sp(2m)}
\labelLL(0,.5)[.3]{m\ge2}
\labelO(1,2)[0]{\sp(2m{+}2)}
\labelU(2,1)[-.4]{\sl(2)}
\labelRR(4,1.5)[0]{\sl(n)}
\labelR(4,1)[0]{\sp(2n)}
\labelRR(4,.5)[0]{n\ge1}
\end{texdraw}
\hskip20pt plus 1fil
\begin{texdraw}
\strichRU(0,1)
\strichRO(1,0)
\strichRO(2,1)
\strichRU(3,2)
\strichRU(4,1)
\strichRO(5,0)
\kreis(0,1)
\punkt(1,0)
\punkt(5,0)
\kreis(6,1)
\labelLL(0,1.5)[0]{\sl(m)}
\labelL(0,1)[0]{\sp(2m)}
\labelLL(0,.5)[0]{m\ge1}
\labelO(2,1)[-0.4]{\sl(2)}
\labelO(3,2)[0]{\sp(4)}
\labelO(4,1)[0.4]{\sl(2)}
\labelRR(6,1.5)[0]{\sl(n)}
\labelR(6,1)[0]{\sp(2n)}
\labelRR(6,.5)[0]{n\ge1}
\end{texdraw}
\vskip5pt
\hrule
\end{table}

\begin{table}[htbp]
\caption{Inference Rules}
\label{rules}
\hrule
\vskip5pt

\leavevmode
\bdiadraw
\strichU(0,1)
\labelR(0,0)[0]{V}
\kreis(0,1)
\labelO(0,1)[0]{\sl(2)}
\ediadraw
\hskip0pt\hskip10pt\raise10pt\hbox{$\Longrightarrow$}\hskip10pt
\bdiadraw
\strichU(0,1)
\strichRO(0,1)
\strichRU(1,2)
\labelR(0,0)[0]{V}
\labelL(0,1)[0]{\sl(2)}
\labelO(1,2)[0]{\sp(2m+2)}
\labelR(2,1)[0]{\sp(2m)}
\labelRR(2,.6)[0]{m\ge1}
\ediadraw
\hskip50pt
\bdiadraw
\strichU(0,1)
\kreis(0,1)
\labelO(0,1)[0]{\sp(4)}
\labelR(0,0)[0]{\omega_1}
\ediadraw
\hskip0pt\hskip10pt\raise10pt\hbox{$\Longrightarrow$}\hskip10pt
\bdiadraw
\strichU(0,1)
\strichU(0,2)
\labelR(0,1)[0]{\sl(4)}
\labelRR(0,0)[.2]{\omega_1}
\labelR(0,0)[-.2]{\omega_3}
\labelR(0,2)[0]{\so(7)}
\ediadraw
\hskip0pt\hskip10pt\raise10pt\hbox{$\isom$}\hskip10pt
\bdiadraw
\strichU(0,1)
\strichU(0,2)
\labelR(0,1)[0]{\so(6)}
\labelRR(0,0)[.2]{\omega_2}
\labelR(0,0)[-.2]{\omega_3}
\labelR(0,2)[0]{\so(7)}
\ediadraw

\vskip5pt
\hrule
\end{table}

\newpage

\begin{table}[htbp]
  \caption{Table of all primitive spherical pairs $(\fg,\fh)$ with
  $\fg$ simple. A ``$*$'' means that $\fs$ is non-trivial in which case
  $\fs=\ft^1$. The third column lists the principal subalgebra $\fl$ of
  the pair and in some cases an indication of its embedding into
  $\fh$. The last column contains marks for reference in the proof of
  Proposition~\ref{ActualComp}.}
\label{KraemerTab}
\baselineskip14pt
{\thickmuskip\medmuskip\medmuskip\thinmuskip\thinmuskip0mu
\halign{$#$\hfil\quad\quad&$#$\hfil\quad\quad&$#$\hfil\quad\quad
&$#$\hfil\quad\quad&$#$\hfil\cr
\fg&\fh&\fl\cr
\noalign{\smallskip\hrule\smallskip}
\sl(m+n)&\*\sl(m)\times\sl(n)&\gl(m-n)\times\ft^{n-1}&m\ge n\ge1&\bullet\cr
\sl(2n)&\sp(2n)&\sl(2)^n&n\ge 2&\bullet\cr
\sl(2n+1)&\*\sp(2n)&\ft^1&n\ge2&b\cr
\sl(n)&\so(n)&0&n\ge3&b\cr
\noalign{\smallskip\hrule\smallskip}
\sp(2m+2n)&\sp(2m)\times\sp(2n)&\sp(2m-2n)\times\sl(2)^n&m\ge n\ge1&\bullet\cr
\sp(2n+2)&\*\sp(2n)&\sp(2n-2)&n\ge1&+\cr
\sp(2n)&\*\sl(n)&0&n\ge2&b\cr
\noalign{\smallskip\hrule\smallskip}
\so(m+n)&\so(m)\times\so(n)^\dagger&\so(m-n)&\llap{\Big\{}
    \vcenter{\ialign{$#$\hfill\cr m\ge n\ge1\cr m+n\ge7\cr}}&\circ\cr
\so(2n)&\*\sl(n)&\llap{\Big\{}\vcenter{\ialign{$#$\hfill&;\ $#$\hfill\cr
    \sl(2)^m&n=2m\cr\sl(2)^m\times\ft^1&n=2m+1\cr}}&n\ge4&\bullet\cr
\so(2n+1)&\*\sl(n)&0&n\ge3&b\cr
\so(7)&\ssG_2&\sl(3)&&b\cr
\so(8)&\ssG_2&\sl(2)&&b\cr
\so(9)&\so(7)^\ddagger&\sl(3)&&b\cr
\so(10)&\*\so(7)^\ddagger&\sl(2)&&b\cr
\noalign{\smallskip\hrule\smallskip}
\ssG_2&\ssA_2&\ssA_1&&\bullet\cr
\ssG_2&\ssA_1\times\ssA_1&0&&b\cr
\ssF_4&\ssB_4&\ssB_3&&x\cr
\ssF_4&\ssC_3\times\ssA_1&0&&b\cr
\ssE_6&\ssC_4&0&&b\cr
\ssE_6&\ssF_4&\ssD_4&&x\cr
\ssE_6&\*\ssD_5&\ssA_3\times\ft^1&&x\cr
\ssE_6&\ssA_5\times\ssA_1&\ft^2\inn\ssA_5&&b'\cr
\ssE_7&\*\ssE_6&\ssD_4&&b\cr
\ssE_7&\ssA_7&0&&b\cr
\ssE_7&\ssD_6\times\ssA_1&\ssA_1^3\inn\ssD_6&&b'\cr
\ssE_8&\ssD_8&0&&b\cr
\ssE_8&\ssE_7\times\ssA_1&\ssD_4\inn\ssE_7&&b'\cr
\noalign{\smallskip\hrule\smallskip}}}

\noindent
${}^\dagger$Read $*\so(m)$ if $n=2$.\hfill\break
\noindent
${}^\ddagger$Embedding via $\so(7)=
\mf{spin}(7)\into\so(8)\into\so(8+\epsilon)$.\hfill\break
\end{table}

\newpage

\begin{table}[htbp] 
  \caption{Table of all primitive spherical pairs $(\fg,\fh)$ with
  $\fg$ not simple. The ``$*$'' indicates $\fs\ne0$ and its embedding
  into $\fg$. The second column lists the principal subalgebra
  $\fl$ of the pair. Again, the last column contains marks for reference in
  the proof of Lemma~\ref{ActualComp}.}
\label{BrionTab}
{\lineskip5pt\thickmuskip\medmuskip\medmuskip\thinmuskip\thinmuskip0mu
\halign{$#$\quad\hfill&$#$\quad\hfill&$#$\hfill\cr
(\fg,\ \fh)&\fl&\cr
\noalign{\smallskip\hrule\smallskip}
\bdiadraw 
\strichRU(0,1)
\strichRO(1,0)
\labelO(0,1)[0]{\fg}
\labelO(2,1)[0]{\fg}
\labelU(1,0)[0]{\fg}
\ediadraw
\scriptstyle\fg\hbox{ \small{simple}}
   &\hbox{max.~torus}&\circ\cr
\bdiadraw
\strichRU(0,1)
\strichRO(1,0)
\Lstern(1,0)
\labelO(0,1)[0]{\sl(n{+}1)}
\labelO(2,1)[0]{\sl(n)}
\labelU(1,0)[0]{\sl(n)}
\ediadraw\scriptstyle n\ge2
   &0&b\cr
\bdiadraw
\strichRU(0,1)
\strichRO(1,0)
\labelO(0,1)[0]{\so(n{+}1)}
\labelO(2,1)[0]{\so(n)}
\labelU(1,0)[0]{\so(n)}
\ediadraw\scriptstyle n\ge5
   &0&b\cr
\bdiadraw
\strichRO(0,0)
\strichRU(1,1)
\strichRO(2,0)
\labelU(0,0)[0]{\sp(2n)}
\labelO(1,1)[0]{\sp(2n{+}4)}
\labelU(2,0)[0]{\sp(4)}
\labelO(3,1)[0]{\sp(4)}
\ediadraw\scriptstyle n\ge1
   &\sp(2n-4)&b'\cr
\bdiadraw
\strichRO(0,0)
\strichRU(1,1)
\strichRO(2,0)
\strichRU(3,1)
\labelU(0,0)[0]{\sl(m)}
\labelO(1,1)[-.1]{\sl(m{+}2)}
\labelU(2,0)[0]{\sl(2)}
\labelO(3,1)[.1]{\sp(2n{+}2)}
\labelU(4,0)[0]{\sp(2n)}
\Mstern(1,1)
\ediadraw\lijst{m\ge1\cr n\ge0\cr}
   &\gl(m-2)\times\sp(2n-2)&x\cr
\bdiadraw
\strichRO(0,0)
\strichRU(1,1)
\strichRO(2,0)
\strichRU(3,1)
\labelU(0,0)[0]{\sp(2m)}
\labelO(1,1)[-.2]{\sp(2m{+}2)}
\labelU(2,0)[0]{\sl(2)}
\labelO(3,1)[.2]{\sp(2n{+}2)}
\labelU(4,0)[0]{\sp(2n)}
\ediadraw\lijst{m\ge0\cr n\ge0\cr}
   &\sp(2m-2)\times\ft^1\times\sp(2n-2)&\bullet\cr
\bdiadraw
\strichRO(0,0)
\strichRU(1,1)
\strichRO(2,0)
\strichRU(3,1)
\strichRO(4,0)
\strichRU(5,1)
\labelU(0,0)[0]{\sp(2m)}
\labelO(1,1)[0]{\sp(2m{+}2)}
\labelU(2,0)[0]{\sl(2)}
\labelO(3,1)[0]{\sp(4)}
\labelU(4,0)[0]{\sl(2)}
\labelO(5,1)[0]{\sp(2n{+}2)}
\labelU(6,0)[0]{\sp(2n)}
\ediadraw\lijst{m\ge0\cr n\ge0\cr}
   &\sp(2m-2)\times\sp(2n-2)&b'\cr
\bdiadraw
\strichRO(0,0)
\move(1 1)\rlvec(4 -1)
\strichRO(2,0)
\move(3 1)\rlvec(2 -1)
\strichRO(4,0)
\strichU(5,1)
\labelU(0,0)[0]{\sp(2l)}
\labelO(1,1)[-.2]{\sp(2l{+}2)}
\labelU(2,0)[0]{\sp(2m)}
\labelO(3,1)[0]{\sp(2m{+}2)}
\labelU(4,0)[-.2]{\sp(2n)}
\labelO(5,1)[.2]{\sp(2n{+}2)}
\labelU(5,0)[.2]{\sl(2)}
\ediadraw\lijst{l\ge0\cr m\ge0\cr n\ge0\cr}
   &\sp(2l-2)\times\sp(2m-2)\times\sp(2n-2)&b'\cr
\noalign{\smallskip\hrule\smallskip}}}
\end{table}

\newpage

\begin{table}[htbp] 
\caption{Table of base components of primitive spherical triples with
  $V\ne0$ and their principal subalgebras. A ``$*$'' indicates
  $\fs\ne0$.}
\label{PairsForTriples}
\hrule
\vskip5pt

\noindent
\bdiadraw
\strichRO(0,1)
\strichRU(1,2)
\strichU(2,1)
\strichU(0,1)
\move(0 1) \rlvec(1.91 -.96)
\KlabelU(2,0)[0]{\ft^{n-1}}
\labelL(0,1)[0]{\sl(m)}
\labelO(1,2)[0]{\sl(m{+}n)}
\labelR(2,1)[0]{\sl(n)}
\KlabelU(0,0)[0]{\gl(m-n)}
\Mstern(1,2)
\labelcond(1,-.6)[0]{m\ge n \ge 1}
\ediadraw%
\hskip0pt plus 1fil
\bdiadraw
\strichU(0,2)
\strichU(0,1)
\labelO(0,2)[0]{\sl(2n)}
\labelR(0,1)[0]{\sp(2n)}
\KlabelU(0,0)[0]{\sl(2)^n}
\labelcond(0,-.6)[0]{n \ge 2}
\ediadraw
\hskip0pt plus 1fil
\bdiadraw
\strichRO(0,1)
\strichRU(1,2)
\strichU(2,1)
\strichU(0,1)
\move(0 1) \rlvec(1.91 -.96)
\KlabelU(2,0)[0]{\sl(2)^{n}}
\labelL(0,1)[0]{\sp(2m)}
\labelO(1,2)[0]{\sp(2m{+}2n)}
\labelR(2,1)[0]{\sp(2n)}
\KlabelU(0,0)[0]{\sp(2m-2n)}
\labelcond(1,-.6)[0]{m \ge n \ge 1}
\ediadraw
\hskip0pt plus 1fil
\bdiadraw
\strichU(0,2)
\strichU(0,1)
\labelO(0,2)[0]{\so(n+1)}
\labelR(0,1)[0]{\so(n)}
\KlabelU(0,0)[0]{\so(n-1)}
\labelcond(0,-.6)[0]{n\ge6}
\ediadraw
\hskip0pt plus 1fil
\bdiadraw
\strichU(0,2)
\strichU(0,1)
\labelO(0,2)[0]{\so(n+2)}
\labelR(0,1)[0]{\so(n)}
\KlabelU(0,0)[0]{\so(n-2)}
\labelcond(0,-.6)[0]{n\ge5}
\Vstern(0,1)
\ediadraw
\hskip0pt plus 1fil

\par
\vskip5pt

\noindent
\bdiadraw
\strichU(0,2)
\strichU(0,1)
\strichRU(0,1)
\labelO(0,2)[0]{\so(2n)}
\labelR(0,1)[.1]{\sl(n)}
\KlabelU(0,0)[0]{\sl(2)^m}
\KlabelR(1,0)[0]{\ft^1\, {\rm if}\ n=2m+1}
\labelRR(1,-.4)[0]{0 \, {\rm if}\ n=2m}
\Vstern(0,1)
\labelcond(.5,-.6)[0]{n \ge 4}
\ediadraw
\hskip0pt plus 1fil
\bdiadraw
\strichU(0,2)
\strichU(0,1)
\labelO(0,2)[0]{\ssG_2}
\labelR(0,1)[0]{\sl(3)}
\KlabelU(0,0)[0]{\sl(2)}
\ediadraw
\hskip0pt plus 1fil
\bdiadraw
\strichRU(0,2)
\strichRO(1,1)
\strichU(1,1)
\labelO(0,2)[0]{\sl(n)}
\labelO(2,2)[0]{\sl(n)}
\labelR(1,1)[-.1]{\sl(n)}
\KlabelU(1,0)[0]{\ft^{n-1}}
\labelcond(1,-.6)[0]{n \ge 3}
\ediadraw
\hskip0pt plus 1fil
\bdiadraw
\strichRO(0,1)
\strichRU(1,2)
\strichRO(2,1)
\strichRU(3,2)
\strichU(2,1)
\strichU(4,1)
\strichU(0,1)
\labelL(0,1)[0]{\sp(2m)}
\labelLL(0,.6)[0]{m\ge0}
\labelO(1,2)[-.1]{\sp(2m{+}2)}
\labelR(2,1)[0]{\sl(2)}
\labelO(3,2)[0.1]{\sp(2n{+}2)}
\labelR(4,1)[0]{\sp(2n)}
\labelRR(4,.5)[0]{n\ge0}
\KlabelU(0,0)[0]{\sp(2m-2)}
\KlabelU(2,0)[0]{\ft^1}
\KlabelU(4,0)[0]{\sp(2n-2)}
\move(0 1) \rlvec(1.91 -.96)
\move(2.07 .07) \rlvec(1.93 .93)
\ediadraw
\vskip5pt
\hrule
\end{table}


\begin{thebibliography}{9999}

\bibitem[BR]{BeRa} C.~Benson and G.~Ratcliff, A classification
of multiplicity free actions, {\em J.\ Algebra}, {\bf 181} (1996), 152--186.

\bibitem[Bou]{Bou} N.~Bourbaki, {\em Groupes et alg\`ebres de Lie,
Chapitres 4, 5 et 6}, Hermann, Paris (1975).

\bibitem[Br1]{Br0}M.~Brion, Repr\'esentations exceptionnelles des
  groupes semi-simples, {\em Ann. Sci. \'Ecole Norm. Sup. (4)} {\bf
    18} (1985), 345--387.

\bibitem[Br2]{Brion}M. Brion, Sur l'image de l'application moment,
  {\em S\'eminaire d'alg\`ebre Paul Dubreil et Marie-Paule Malliavin
  (Paris, 1986)}, Lecture Notes in Math., {\bf 1296}, Springer,
  Berlin (1987), 177--192.

\bibitem[Br3]{Br} M.~Brion, Classification des espaces homog\`enes
sph\'eriques, {\em Compos.\ Math.}, {\bf 63} (1987), 189--208.

\bibitem[Ca]{Cam} R.~Camus, Vari\'et\'es sph\'eriques affines lisses,
  {\em Th\`ese de doctorat (Universit\'e J. Fourier)} (2001)
{\tt www-fourier.ujf-grenoble.fr/THESE/html/a110}

\bibitem[De]{Del} Th.~Delzant, Classification des actions
hamiltoniennes compl\`etement int\'egrables de rang deux, {\em Ann.
Global Anal. Geom.}  {\bf 8} (1990), 87--112.

\bibitem[El1]{ElasCan} A.G.~\'Elashvili, Canonical form and stationary
subalgebras of points of general position for simple linear Lie groups
(Russian), {\em Funkcional.\ Anal.\ i Prilo\v zen.}, {\bf 6} (1972),
51--63. English translation: {\em Functional Anal.\ Appl.}, {\bf 6}
(1972), 44--53.

\bibitem[El2]{ElasStat} A.G.~\'Elashvili, Stationary subalgebras of
points of the common state for irreducible linear Lie groups
(Russian), {\em Funkcional.\ Anal.\ i Prilo\v zen.}, {\bf 6} (1972),
65--78. English translation: {\em Functional Anal.\ Appl.}, {\bf 6}
(1972), 139--148.

\bibitem[Kac]{Kac}V.~Kac, Some
remarks on nilpotent orbits, {\em J.~Algebra}, {\bf 64} (1980), 190--213.

\bibitem[Kn1]{WM} F.~Knop, Weylgruppe und Momentabbildung, {\em
Invent.\ Math.}, {\bf 99} (1990), 1--23.

\bibitem[Kn2]{IB} F.~Knop, \"Uber Bewertungen, welche unter einer
  reduktiven Gruppe invariant sind, {\em Mathematische Annalen} {\bf
  295} (1993), 333-363.

\bibitem[Kn3]{AsBh}F.~Knop, The asymptotic behavior of invariant
collective motion, {\em Invent.\ Math.}, {\bf 116} (1994), 309--328.

\bibitem[Kn4]{KRem}F.~Knop, Some remarks on multiplicity free spaces,
Proc.\ NATO Adv.\ Study Inst.\ on Representation Theory and Algebraic
Geometry (A.~Broer, G.~Sabidussi, eds.), {\em Nato ASI Series C}, {\bf
514}, Kluwer, Dordrecht (1998), 301--317.

\bibitem[Kr]{Kra} M.~Kr\"amer, Sph\"arische Untergruppen in
kompakten zusammenh\"angenden Liegruppen, {\em Compos.\
Math.}, {\bf 38} (1979), 129--153.

\bibitem[Le]{Le} A.~Leahy, A classification of multiplicity
free representations, {\em J.\ Lie Theory}, {\bf 8} (1998), 367--391.

\bibitem[Lu1]{LunaSl} D.~Luna, Slices \'etales, {\em M\`emoires de la
S.M.F.}, {\bf 33} (1973), 81--105.

\bibitem[Lu2]{Luna} D.~Luna, Vari\'et\'es sph\'eriques de type
  $\ssA$, {\em Publ. Math. Inst. Hautes \'Etudes Sci.} {\bf 94}
(2001), 161--226.

\bibitem[Mi]{Mik}I.V.~Mikityuk, On the integrability of
invariant Hamiltonian systems with homogeneous configuration
spaces, {\em Math.\ USSR--Sb.}, {\bf 57} (1987), 527--546.

\bibitem[Pa1]{Pan1}D.I.~Panyushev, 
Complexity and rank of homogeneous spaces,
{\em Geom. Dedicata} {\bf34} (1990), 249--269.

\bibitem[Pa2]{Pan2}D.I.~Panyushev, Complexity and rank of actions
    in invariant theory. Algebraic geometry, 8, {\em J. Math. Sci. (New
  York)}  {\bf95}  (1999), 1925--1985

\bibitem[Ri]{Rich} R.W.~Richardson, Principal orbit types for algebraic
transformation spaces in characteristic zero, {\em Invent.\ Math.},
{\bf 16} (1972), 6--14.
      
\bibitem[Sch]{Schwarz}G.W.~Schwarz, Representations of simple Lie groups
  with regular rings of invariants, {\em Invent.\ Math.}, {\bf 49}
  (1978), 167--191.

\bibitem[Sj]{Sja}R.~Sjamaar, Convexity properties of the moment map
re-examined, {\em Adv. Math.}, {\bf138} (1998), 46--91.

\bibitem[VK]{VinKim}E.B.~Vinberg, B.~Kimelfeld, Homogeneous domains on
  flag manifolds and spherical subgroups of semisimple Lie groups,
  {\em Funktsional. Anal. i Prilozhen.}, {\bf12} (1978), 12--19.
  English translation: {\em Functional Anal. Appl.} {\bf 12} (1978),
    168--174 (1979).

\end{thebibliography}
\end{document}